\documentclass[12pt,reqno]{amsart}
\usepackage{amsmath,amssymb,latexsym, amsthm, amscd, mathrsfs, stmaryrd} 
\usepackage[linktocpage=true]{hyperref}
\hypersetup{colorlinks,linkcolor=blue,urlcolor=cyan,citecolor=blue}

\usepackage[all]{xy}
\usepackage{hyperref}
\usepackage{color}
\newcommand{\nc}{\newcommand}
\nc{\browntext}[1]{\textcolor{brown}{#1}}
\nc{\greentext}[1]{\textcolor{green}{#1}}
\nc{\redtext}[1]{\textcolor{red}{#1}}
\nc{\bluetext}[1]{\textcolor{blue}{#1}}
\nc{\brown}[1]{\browntext{ #1}}
\nc{\green}[1]{\greentext{ #1}}
\nc{\red}[1]{\redtext{ #1}}
\nc{\blue}[1]{\bluetext{ #1}}
\nc{\zb}[1]{\redtext{From zb: #1}}

\newtheorem{theorem}{Theorem}  [section]
\newtheorem{corollary}[theorem]{Corollary}
\newtheorem{lemma}[theorem]{Lemma}

\newtheorem{proposition}[theorem]{Proposition}

\newtheorem{definition}[theorem]{Definition}

\theoremstyle{remark}
\newtheorem{remark}[theorem]{Remark}

\numberwithin{equation}{section}

\def\haTh{\widehat{\Theta}}
\def \haH{\widehat{H}}

\newcommand{\bH}{\mathbf H}

\def \E{K}

\def \hath{\widehat{\theta}}
\def \bt{\mathbf t}
\def \bn{\mathbf n}
\def \bh{\mathbf h}
\newcommand{\cc}{{\mathcal C}}
\def \bC{\mathbf C}
\def \dB{\Theta}
\def \bpi{\boldsymbol{\pi}}

\def \co{\mathcal O}

\def \ch{{\mathcal H}}

\def \ct{{\mathcal T}}

\def \haB{\widehat{B}}
\newcommand{\tMHX}{{}^\imath\widetilde{\ch}(\X_\bfk)}
\newcommand{\tCMH}{{}^\imath\widetilde{\cc}(\bfk Q)}
\newcommand{\tCMHX}{{}^\imath\widetilde{\cc}(\X_\bfk)}
\renewcommand{\mod}{\operatorname{mod}\nolimits}
\newcommand{\tCMHC}{{}^\imath\widetilde{\cc}(\bfk C_n)}
\newcommand{\tCHX}{\widetilde{\cc}(\X_\bfk)}

\numberwithin{equation}{section}

\renewcommand{\ker}{\operatorname{Ker}\nolimits}

\def \cp{\mathcal P}
\newcommand{\Hom}{\operatorname{Hom}\nolimits}

\newcommand{\Aut}{\operatorname{Aut}\nolimits}
\newcommand{\Id}{\operatorname{Id}\nolimits}

\newcommand{\Ext}{\operatorname{Ext}\nolimits}

\def \y{{B}}
\def \haB{\widehat{B}}

\newcommand{\mbf}{\mathbf}

\newcommand{\mrm}{\mathrm}

\newcommand{\End}{\mrm{End}}

\newcommand{\de}{\delta}

\newcommand{\N}{\mathbb N}

\newcommand{\qbinom}[2]{\begin{bmatrix} #1\\#2 \end{bmatrix} }
\newcommand{\Q}{\mathbb Q}
\newcommand{\sll}{\mathfrak{sl}}
\newcommand{\T}{\texttt{\rm T}}

\newcommand{\U}{\mbf U}
\newcommand{\K}{\mathbb K}
\newcommand{\F}{\mathbb F}

\newcommand{\bs}{\mathbf s}
\newcommand{\arxiv}[1]{\href{http://arxiv.org/abs/#1}{\tt arXiv:\nolinkurl{#1}}}

\newcommand{\Z}{\mathbb Z}

\newcommand{\TT}{\mathbf T}

\def \X{\mathbb X}
\def \fg{\mathfrak{g}}
\def \bU{{\mathbf U}}
\newcommand{\tK}{\widetilde{K}}
\def \I{\mathbb{I}}
\def \bv{v}

\newcommand{\tUiD}{{}^{\text{Dr}}\tUi}

\def \hn{\widehat{\mathfrak{n}}}

\newcommand{\tUi}{\widetilde{{\mathbf U}}^\imath}
\newcommand{\sqq}{{\bf v}}

\newcommand{\dbl}{\operatorname{dbl}\nolimits}
\newcommand{\swa}{\operatorname{swap}\nolimits}

\newcommand{\coker}{\operatorname{Coker}\nolimits}

\newcommand{\tU}{\widetilde{\mathbf U}}

\def \cl{L}

\def \K{\mathbb{K}}
\def \R{\mathbb{R}}
\def \SS{\mathbb{S}}

\def \cc{\mathcal C}
\def\ca{\mathcal A}
\newcommand{\Iso}{\operatorname{Iso}\nolimits}
\renewcommand{\Im}{\operatorname{Im}\nolimits}
\newcommand{\res}{\operatorname{res}\nolimits}
\newcommand{\iH}{{}^\imath\widetilde{\ch}}

\newcommand{\coh}{\operatorname{coh}\nolimits}
\newcommand{\rep}{\operatorname{rep}\nolimits}

\def \PL{\mathbb{P}^1_{\bfk}}

\def \scrf{\mathscr F}
\def \scrt{\mathscr T}
\def \P{\mathbb P}
\def \cI{\mathcal I}

\def \cs{{\mathcal{S}}}
\def \ch{\mathcal H}
\def \cd{\mathcal D}
\def\bfk{\mathbf{k}}

\def \bp{\mathbf p}
\def \ul{\underline}

\def \II{\I_0}
\def \Lg{L\fg}

\newcommand{\haT}{\widehat{\Theta}}

\def \bH{\mathbf{H}}
\newcommand{\gr}{\operatorname{gr}\nolimits}

\def \cR{\mathcal R}
\def \ce{\mathcal E}

\def \bla{\boldsymbol{\lambda}}
\def \blx{x}

\allowdisplaybreaks

\begin{document}

	\title[$\imath$Hall algebras of weighted projective lines and quantum symmetric pairs II]
	{$\imath$Hall algebras of weighted projective lines and quantum symmetric pairs II: injectivity}
	
	\author[Ming Lu]{Ming Lu}
	\address{Department of Mathematics, Sichuan University, Chengdu 610064, P.R.China}
	\email{luming@scu.edu.cn}

	\author[Shiquan Ruan]{Shiquan Ruan}
	\address{ School of Mathematical Sciences,
		Xiamen University, Xiamen 361005, P.R.China}
	\email{sqruan@xmu.edu.cn}


	\subjclass[2020]{Primary 17B37, 16E60, 18E35.}
	
	\keywords{Quantum symmetric pairs, Affine $\imath$quantum groups, Drinfeld type presentations, Hall algebras, Weighted projective lines}
	
	\begin{abstract}
		We show that the morphism $\Omega$ 
		from the $\imath$quantum loop algebra $\tUiD(\Lg)$ of split type to the $\imath$Hall algebra of the weighted projective line is injective if $\fg$ is of finite or affine type. As a byproduct, we use the whole $\imath$Hall algebra of the cyclic quiver $C_n$ to realise the $\imath$quantum loop algebra of affine $\mathfrak{gl}_n$. 
	\end{abstract}
	
	\maketitle
	\setcounter{tocdepth}{1}
	\tableofcontents

	
	\section{Introduction}
	In 1990s, Ringel \cite{Rin90} used the Hall algebra of Dynkin quiver to realize a half part of a quantum group, which is generalized to Kac-Moody setting by Green \cite{Gr95}. Two decades later, Bridgeland \cite{Br13} introduced a Hall algebra of the category of 2-periodic complexes of an acyclic quiver to realise the entire quantum group.  Inspired by Bridgeland's construction, the first author and Peng \cite{LP16}  defined a so-called semi-derived Ringel-Hall algebra for arbitrary hereditary abelian categories via 2-periodic complexes. The construction of semi-derived Ringel-Hall algebras was further extended to $m$-periodic complexes over heredtiary abelian categories, and then to 1-Gorenstein algebras; see \cite{LP19,LW19a}.

	There has been a current realization of the affine quantum groups formulated by Drinfeld \cite{Dr88, Be94, Da15}. Hall algebra of the projective line was studied in a visionary paper by Kapranov \cite{Ka97} and then extended by Baumann-Kassel \cite{BKa01} to realize the current half of quantum affine $\sll_2$. The Hall algebra of a weighted projective line was developed in \cite{Sch04} to realize half an affine quantum group of ADE type, which were then upgraded to the whole quantum group via Drinfeld double techniques  \cite{DJX12, BS13}.
	
	The universal $\imath$quantum group is introduced in \cite{LW19a}, whose central reduction reproduces the $\imath$quantum group by G. Letzter \cite{Let99}. Recently, a Drinfeld type presentation of  $\imath$quantum groups of split affine ADE type is established in \cite{LW20b}.

	In recent years, Wang and the first author \cite{LW19a, LW20a} have developed $\imath$Hall algebras of $\imath$quivers to realize the (universal) quasi-split $\imath$quantum groups of Kac-Moody type. 
	The $\imath$Hall algebras are constructed in the framework of semi-derived Ringel-Hall algebras of 1-Gorenstein algebras.
	Together with Wang, we \cite{LRW20a} used the $\imath$Hall algebra of the projective line to give a geometric realization of the $q$-Onsager algebra (i.e., $\imath$quantum group of split affine $\mathfrak{sl}_2$) in its Drinfeld type presentation. 
	
	$\imath$Quantum loop algebras  are a further generalization of $\imath$quantum groups of split affine ADE type considered in \cite{LW20b}. In \cite{LR21}, we established a homomorphism $\Omega$ from the $\imath$quantum loop algebra $\tUiD(\Lg)$ of split type to the $\imath$Hall algebra of the weighted projective line $\iH(\X_\bfk)$.

	This is a sequel of \cite{LR21} which is devoted to proving the morphism $\Omega:\tUiD(\Lg)\rightarrow \iH(\X_\bfk)$ given in \cite{LR21} is injective if $\fg$ is of finite or affine type. The idea of the proof is to construct the inverse of $\Omega$ by giving a presentation of the composition subalgebra $\tCMHX$ of $\iH(\X_\bfk)$. This presentation of $\tCMHX$ is established by using the presentation of the ``positive half" of the quantum loop algebra $\U_v(\widehat{\bn})$ and the isomorphism between $\U_v(\widehat{\bn})$ and  $\widetilde{\cc}(\X_\bfk)$ given in \cite{Sch04}, since  $\iH(\X_\bfk)$ is naturally a filtered algebra with its associated grading algebra isomorphic to $\widetilde{\ch}(\X_\bfk)$. In this procedure, we need to use the entire $\imath$Hall algebra of the cyclic quiver $C_n$ to realize the $\imath$quantum loop algebra of affine $\mathfrak{gl}_n$; see \cite{Sch02,Hu05} for the corresponding results for quantum loop algebra of affine $\mathfrak{gl}_n$.
	
	This paper is organized as follows. In Section \ref{sec:prel}, we recall the necessary materials on $\imath$quantum groups, Hall algebras and $\imath$Hall algebras, weighted projective lines, etc. In Section \ref{sec:hom}, we recall the construction of the morphism $\Omega$ given in \cite{LR21}, and formulate the main result. Section \ref{sec:Hall cyclic} is devoted to realizing the $\imath$quantum loop algebra of affine $\mathfrak{gl}_n$ via the $\imath$Hall algebra of the cyclic quiver $C_n$. In Section \ref{sec:injectivity}, we give a presentation of the composition subalgebra of $\iH(\X_\bfk)$, and then prove the main result.
	
	\subsection{Acknowledgments}
	ML deeply thanks Weiqiang Wang for guiding him to study the $\imath$quantum groups, and also his continuing encouragement.
	ML thanks University of Virginia for hospitality and support. We thank Liangang Peng and Jie Xiao for helpful discussions on Ringel-Hall algebras of weighted projective lines. ML is partially supported by the National Natural Science Foundation of China (No. 12171333).
	SR is partially supported by the National Natural Science Foundation of China (No. 11801473) and the Fundamental Research Funds for Central Universities of China (No. 20720210006).

	\section{Preliminaries}
	\label{sec:prel}
	The preliminaries are the same as in \cite{LR21},  we recall them briefly to make this paper self-contained.

	\subsection{$\imath$Quantum groups}
	
	For $n\in \Z, r\in \N$, denote by
	\[
	[n] =\frac{v^n -v^{-n}}{v-v^{-1}},\qquad
	\qbinom{n}{r} =\frac{[n][n-1]\ldots [n-r+1]}{[r]!}.
	\]
	For $A, B$ in a $\Q(v)$-algebra, we shall denote $[A,B]_{v^a} =AB -v^aBA$, and $[A,B] =AB - BA$. We also use the convenient notation throughout the paper
	\begin{align}
		[A,B,C]_{v^a}=\big[[A,B]_{v^a},C\big]_{v^a}.
	\end{align}

	Let $C=(c_{ij})_{i,j\in \I}$ be the generalized Cartan matrix (GCM) of a Kac-Moody algebra $\fg$. 
	The notion of (quasi-split) universal $\imath$quantum groups $\tUi$ was formulated in \cite{LW19a}.
	
	Let $C=(c_{ij})_{i,j\in \I}$ be the GCM of an affine Lie algebra $\widehat{\fg}$, for $\I = \II \cup \{0\}$ with the affine node $0$.
	The {\em universal affine $\imath$quantum group of split type ADE} is the $\Q(v)$-algebra $\tUi =\tUi(\widehat{\fg})$ with generators $B_i$, $\K_i^{\pm 1}$ $(i\in \I)$, subject to the following relations, for $i, j\in \I$:
	\begin{eqnarray}
		\label{eq:KK}
		\K_i\K_i^{-1} =\K_i^{-1}\K_i=1, \qquad\qquad \quad \K_i  \text{ is central},&
		\\
		B_iB_j -B_j B_i=0, \qquad\qquad\qquad\qquad\qquad &&\text{ if } c_{ij}=0,
		\label{eq:S1} \\
		B_i^2 B_j -[2] B_i B_j B_i +B_j B_i^2 = - v^{-1}  B_j \K_i,  \qquad\qquad\qquad &&\text{ if }c_{ij}=-1,
		\label{eq:S2} \\
		\sum_{r=0}^3 (-1)^r \qbinom{3}{r} B_i^{3-r} B_j B_i^{r} = -v^{-1} [2]^2 (B_iB_j-B_jB_i) \K_i,  & & \text{ if }c_{ij}=-2.
		\label{eq:S3}
	\end{eqnarray}
	
	Let $\alpha_i$ $(i\in \I)$ be the simple roots of $\widehat{\fg}$, and $\alpha_0=\de -\theta$, where $\de$ denotes the basic imaginary root, and $\theta$ is the highest root of $\fg$. 	For $\mu = \mu' +\mu''  \in \Z \I := \oplus_{i\in \I} \Z \alpha_i$,  define $\K_\mu$ such that
	\begin{align}
		\K_{\alpha_i} =\K_i, \quad
		\K_{-\alpha_i} =\K_i^{-1}, \quad
		\K_{\mu} =\K_{\mu'} \K_{\mu''},
		\quad  \K_\delta =\K_0 \K_\theta.
	\end{align}

	\begin{lemma} [\cite{LW21b}; \text{also cf. \cite{KP11, BK20}}]
		\label{lem:Ti}
		For $i\in \I$, there exists an automorphism $\TT_i$ of the $\Q(v)$-algebra $\tUi$ such that
		$\TT_i(\K_\mu) =\K_{s_i\mu}$ for $\mu\in \Z\I$, and
		\[
		\TT_i(B_j)= \begin{cases}
			\K_i^{-1} B_i,  &\text{ if }j=i,\\
			B_j,  &\text{ if } c_{ij}=0, \\
			B_jB_i-vB_iB_j,  & \text{ if }c_{ij}=-1, \\
			{[}2]^{-1} \big(B_jB_i^{2} -v[2] B_i B_jB_i +v^2 B_i^{2} B_j \big) + B_j\K_i,  & \text{ if }c_{ij}=-2,
		\end{cases}
		\]
		for $j\in \I$.
		Moreover,  $\TT_i$ $(i\in \I)$ satisfy the braid group relations, i.e., $\TT_i \TT_j =\TT_j \TT_i$ if $c_{ij}=0$, and $\TT_i \TT_j \TT_i =\TT_j \TT_i \TT_j$ if $c_{ij}=-1$.
	\end{lemma}

	
	
	\subsection{A Drinfeld type presentation}
	\label{subsec:Dr2}
	
	Let $C=(c_{ij})_{i,j\in \II}$ be a GCM of a {\em simply-laced} Kac-Moody algebra $\fg$. The $\imath$quantum loop algebra $\tUiD$ of split type is the $\Q(v)$-algebra  generated by $\K_{i}^{\pm1}$, $C^{\pm1}$, $H_{i,m}$, $\Theta_{i,l}$ and $\y_{i,l}$, where  $i\in \II$, $m \in \Z_{+}$, $l\in\Z$, subject to some relations.
	In order to give the explicit definition of $\tUiD$,  we introduce some shorthand notations below.
	
	Let $k_1, k_2, l\in \Z$ and $i,j \in \II$. Set
	\begin{align}
		\begin{split}
			S(k_1,k_2|l;i,j)
			&=  B_{i,k_1} B_{i,k_2} B_{j,l} -[2] B_{i,k_1} B_{j,l} B_{i,k_2} + B_{j,l} B_{i,k_1} B_{i,k_2},
			\\
			\SS(k_1,k_2|l;i,j)
			&= S(k_1,k_2|l;i,j)  + \{k_1 \leftrightarrow k_2 \}.
			\label{eq:Skk}
		\end{split}
	\end{align}
	Here and below, $\{k_1 \leftrightarrow k_2 \}$ stands for repeating the previous summand with $k_1, k_2$ switched if $k_1\neq k_2$, so the sums over $k_1, k_2$ are symmetric.
	We also denote
	\begin{align}
		\begin{split}
			R(k_1,k_2|l; i,j)
			&=   \K_i  C^{k_1}
			\Big(-\sum_{p\geq0} v^{2p}  [2] [\Theta _{i,k_2-k_1-2p-1},\y_{j,l-1}]_{v^{-2}}C^{p+1}
			\label{eq:Rkk} \\
			&\qquad\qquad -\sum_{p\geq 1} v^{2p-1}  [2] [\y_{j,l},\Theta _{i,k_2-k_1-2p}]_{v^{-2}} C^{p}
			- [\y_{j,l}, \Theta _{i,k_2-k_1}]_{v^{-2}} \Big),
			\\
			\R(k_1,k_2|l; i,j) &= R(k_1,k_2|l;i,j) + \{k_1 \leftrightarrow k_2\}.
		\end{split}
	\end{align}
	Sometimes, it is convenient to rewrite part of the summands in \eqref{eq:Rkk} as
	\begin{align*}
		&-\sum_{p\geq 1} v^{2p-1}  [2] [\y_{j,l},\Theta _{i,k_2-k_1-2p}]_{v^{-2}} C^{p}
		- [\y_{j,l}, \Theta _{i,k_2-k_1}]_{v^{-2}}\\
		=&
		-\sum_{p\geq 0} v^{2p-1}  [2] [\y_{j,l},\Theta _{i,k_2-k_1-2p}]_{v^{-2}} C^{p}
		+v^{-2}[\y_{j,l}, \Theta _{i,k_2-k_1}]_{v^{-2}}.
	\end{align*}
	
	
	\begin{definition}[$\imath$quantum loop algebras, \cite{LW20b}]
		\label{def:iDR}
		Let $C=(c_{ij})_{i,j\in \II}$ be a GCM of a simply-laced Kac-Moody algebra $\fg$. 
		The $\imath$quantum loop algebra $\tUiD$ of split type is the $\Q(v)$-algebra  generated by $\K_{i}^{\pm1}$, $C^{\pm1}$, $H_{i,m}$ and $\y_{i,l}$, where  $i\in \II$, $m \in \Z_{+}$, $l\in\Z$, subject to the following relations (for $m,n \in \Z_{+}$ and $k,l\in \Z$):
		\begin{align}
			& \K_i, C \text{ are central, } \quad \K_i\K_i^{-1}=1, \;\; C C^{-1}=1,
			\label{iDR1a}
			\\
			&[H_{i,m},H_{j,n}]=0,\label{iDR1b}\\
			&[H_{i,m},\y_{j,l}]=\frac{[mc_{ij}]}{m} \y_{j,l+m}-\frac{[mc_{ij}]}{m} \y_{j,l-m}C^m,
			\label{iDR2}
			\\
			&[\y_{i,k} ,\y_{j,l}]=0,   \text{ if }c_{ij}=0,  \label{iDR4}
			\\
			&[\y_{i,k}, \y_{j,l+1}]_{v^{-c_{ij}}}  -v^{-c_{ij}} [\y_{i,k+1}, \y_{j,l}]_{v^{c_{ij}}}=0, \text{ if }i\neq j,
			\label{iDR3a}
			\\ 
			&[\y_{i,k}, \y_{i,l+1}]_{v^{-2}}  -v^{-2} [\y_{i,k+1}, \y_{i,l}]_{v^{2}}
			=v^{-2}\Theta_{i,l-k+1} C^k \K_i-v^{-4}\Theta_{i,l-k-1} C^{k+1} \K_i
			\label{iDR3b} \\
			&\qquad\qquad\qquad\qquad\qquad\qquad\quad\quad\quad
			+v^{-2}\Theta_{i,k-l+1} C^l \K_i-v^{-4}\Theta_{i,k-l-1} C^{l+1} \K_i, \notag
			\\
			\label{iDR5}
			&    \SS(k_1,k_2|l; i,j) =  \R(k_1,k_2|l; i,j), \text{ if }c_{ij}=-1.
		\end{align}
		Here we set
		\begin{align}  \label{Hm0}
			{\Theta}_{i,0} =(v-v^{-1})^{-1}, \qquad {\Theta}_{i,m} =0, \; \text{ for }m<0;
		\end{align}
		and $\Theta_{i,m}$ ($m\geq1$)  are related to  $H_{i,m}$ by the following equation:
		\begin{align}
			\label{exp h}
			1+ \sum_{m\geq 1} (v-v^{-1})\Theta_{i,m} u^m  = \exp\Big( (v-v^{-1}) \sum_{m\geq 1} H_{i,m} u^m \Big).
		\end{align}
	\end{definition}
	Let us mention that in spite of its appearance, the RHS of \eqref{iDR3b} typically has two nonzero terms, thanks to the convention \eqref{Hm0}.

	Let $C$ be a GCM of an affine Lie algebra $\fg$ of type ADE.
	Define a sign function
	\begin{align}
		\label{eq:sign}
		o(\cdot): \,\II \longrightarrow \{\pm 1\}
	\end{align}
	such that $o(i) o(j)=-1$ whenever $c_{ij} <0$ (there are clearly exactly 2 such functions).
	We define the following elements in $\tUi$, for $i\in \II$, $k\in \Z$ and $m\ge 1$: 
	\begin{align}
		B_{i,k} &= o(i)^k \TT_{\omega_i}^{-k} (B_i),
		\label{Bik} \\
		\acute{\Theta}_{i,m} &=  o(i)^m \Big(-B_{i,m-1} \TT_{\omega_i'} (B_i) +v^{2} \TT_{\omega_i'} (B_i) B_{i,m-1}
		\label{Thim1} \\
		& \qquad\qquad\qquad\qquad + (v^{2}-1)\sum_{p=0}^{m-2} B_{i,p} B_{i,m-p-2}  \K_{i}^{-1}\K_{\de} \Big),
		\notag \\
		\dB_{i,m} &=\acute{\Theta}_{i,m} - \sum\limits_{a=1}^{\lfloor\frac{m-1}{2}\rfloor}(v^2-1) v^{-2a} \acute{\Theta}_{i,m-2a}\K_{a\de} -\de_{m,ev}v^{1-m} \K_{\frac{m}{2}\de}.
		\label{Thim}
	\end{align}
	and  $H_{i,m}$ is defined via \eqref{exp h}.
	In particular, $B_{i,0}=B_i$. 
	
	\begin{theorem}[\cite{LW20b}]
		\label{thm:ADE}
		If $C$ is a GCM of an affine Lie algebra $\fg$ of type ADE, then there is a $\Q(v)$-algebra isomorphism ${\Phi}: \tUiD \rightarrow\tUi$, which sends
		\begin{align}   \label{eq:map}
			\y_{i,k}\mapsto B_{i,k},  \quad H_{i,m}\mapsto H_{i,m}, \quad \Theta_{i,m}\mapsto \Theta_{i,m}, \quad
			\K_i\mapsto \K_i, 
			\quad C\mapsto \K_\de,
		\end{align}
		for $i\in \II, k\in \Z$, and $m \ge 1$.
	\end{theorem}
	
	\subsection{Star-shaped graph}
	
	For $\bp=(p_1,\dots,p_\bt)\in\Z_{+}^\bt$, let us consider the following star-shaped graph $\Gamma=\T_{p_1,\dots,p_\bt}$:
	
	\begin{center}\setlength{\unitlength}{0.8mm}
		\begin{equation}
			\label{star-shaped}
			\begin{picture}(110,30)(0,35)
				\put(0,40){\circle*{1.4}}
				\put(2,42){\line(1,1){16}}
				\put(20,60){\circle*{1.4}}
				\put(23,60){\line(1,0){13}}
				\put(40,60){\circle*{1.4}}
				\put(43,60){\line(1,0){13}}
				\put(60,58.5){\large$\cdots$}
				\put(70,60){\line(1,0){13}}
				\put(88,60){\circle*{1.4}}

				\put(3,41){\line(4,1){13}}
				\put(20,45){\circle*{1.4}}
				\put(23,45){\line(1,0){13}}
				\put(40,45){\circle*{1.4}}
				\put(43,45){\line(1,0){13}}
				\put(60,43.5){\large$\cdots$}
				\put(70,45){\line(1,0){13}}
				\put(88,45){\circle*{1.4}}
				
				\put(19,30){\Large$\vdots$}
				
				\put(39,30){\Large$\vdots$}
				
				\put(87,30){\Large$\vdots$}
				
				\put(2,38){\line(1,-1){16}}
				\put(20,20){\circle*{1.4}}
				\put(23,20){\line(1,0){13}}
				\put(40,20){\circle*{1.4}}
				\put(43,20){\line(1,0){13}}
				\put(60,18.5){\large$\cdots$}
				\put(70,20){\line(1,0){13}}
				\put(88,20){\circle*{1.4}}
				
				\put(-4.5,39){$\star$}
				
				\put(16,62){\tiny$[1,1]$}
				
				\put(36,62){\tiny$[1,2]$}
				\put(81,62){\tiny$[1,p_1-1]$}

				\put(16.5,47){\tiny$[2,1]$}
				\put(36.5,47){\tiny$[2,2]$}
				\put(81,47){\tiny$[2,p_2-1]$}

				\put(16.5,16){\tiny$[\bt,1]$}
				\put(36.5,16){\tiny$[\bt,2]$}
				\put(81,16){\tiny$[\bt,p_\bt-1]$}

				
			\end{picture}
		\end{equation}
		\vspace{1cm}
	\end{center}
	The set of vertices is denoted by $\II$. As marked in the graph, the central vertex is denoted by $\star$. Let $J_1,\dots,J_\bt$ be the branches, which are subdiagrams of type $A_{p_1-1},\dots,A_{p_\bt-1}$ respectively. Denote by $[i,j]$ the $j$-th vertex in the $i$-th branch. These examples includes all finite-type Dynkin diagrams as well as the affine Dynkin diagrams of types $D_4^{(1)},E_6^{(1)},E_7^{(1)}$, and $E_8^{(1)}$.
	
	Let $\fg$ be the Kac-Moody algebra corresponding to $\Gamma$, and $\cl \fg$ be its loop algebra. Then we have the $\imath$quantum loop algebra $\tUiD$.
	In view of the graph $\Gamma$, for the root system $\cR$ of $\cl\fg$, its simple roots are denoted by $\alpha_\star$ and $\alpha_{ij}$, for $1\leq i\leq \bt$ and $1\leq j\leq p_i-1$.

	
	\begin{lemma}[\text{\cite[Lemma 2.6]{LR21}}]
		\label{lem:reduced generators}
		For any $a\in\Z$,
		$\tUiD$ is  generated by $\y_{\star,a}$, $\y_{\star,a+1}$, $\y_{i,0}$, $\K_{i}^{\pm1}$ and $C^{\pm1}$ for all  $i\in \II-\{\star\}$.
	\end{lemma}

	

	\subsection{Hall algebras}
	
	Let $\ce$ be an essentially small exact category in the sense of Quillen, linear over  $\bfk=\F_q$.
	Assume that $\ce$ has finite morphism and extension spaces, i.e.,
	\[
	|\Hom(M,N)|<\infty,\quad |\Ext^1(M,N)|<\infty,\,\,\forall M,N\in\ce.
	\]
	
	Given objects $M,N,L\in\ce$, define $\Ext^1(M,N)_L\subseteq \Ext^1(M,N)$ as the subset parameterizing extensions whose middle term is isomorphic to $L$. We define the {\em Ringel-Hall algebra} $\ch(\ce)$ (or {\em Hall algebra} for short) to be the $\Q$-vector space whose basis is formed by the isoclasses $[M]$ of objects $M$ in $\ce$, with the multiplication defined by (see \cite{Br13})
	\begin{align}
		\label{eq:mult}
		[M]\diamond [N]=\sum_{[L]\in \Iso(\ce)}\frac{|\Ext^1(M,N)_L|}{|\Hom(M,N)|}[L].
	\end{align}
	
	For any three objects $L,M,N$, let
	\begin{align}
		\label{eq:Fxyz}
		F_{MN}^L:= \big |\{X\subseteq L \mid X \cong N,  L/X\cong M\} \big |.
	\end{align}
	The Riedtmann-Peng formula states that
	\[
	F_{MN}^L= \frac{|\Ext^1(M,N)_L|}{|\Hom(M,N)|} \cdot \frac{|\Aut(L)|}{|\Aut(M)| |\Aut(N)|}.
	\]

	\subsection{Semi-derived Hall algebras and $\imath$Hall algebras}
	
	\label{subsec:periodic}
	
	Let $\ca$ be a hereditary abelian category  which is essentially small with finite-dimensional homomorphism and extension spaces. 
	
	A $1$-periodic complex $X^\bullet$ in $\ca$ is a pair $(X,d)$ with $X\in\ca$ and a differential $d:X\rightarrow X$. A morphism $(X,d) \rightarrow (Y,e)$ is given by a morphism $f:X\rightarrow Y$ in $\ca$ satisfying $f\circ d=e\circ f$. Let $\cc_1(\ca)$ be the category of all $1$-periodic complexes in $\ca$. Then $\cc_1(\ca)$ is an abelian category.
	A $1$-periodic complex $X^\bullet=(X,d)$ is called acyclic if $\ker d=\Im d$. We denote by $\cc_{1,ac}(\ca)$ the full subcategory of $\cc_1(\ca)$ consisting of acyclic complexes.
	Denote by $H(X^\bullet)\in\ca$ the cohomology group of $X^\bullet$, i.e., $H(X^\bullet)=\ker d/\Im d$, where $d$ is the differential of $X^\bullet$.

	For any $X\in\ca$, denote the stalk complex by
	\[
	C_X =(X,0)
	\]
	(or just by $X$ when there is no confusion), and denote by $K_X$ the following acyclic complex:
	\[
	K_X:=(X\oplus X, d),
	\qquad \text{ where }
	d=\left(\begin{array}{cc} 0&\Id \\ 0&0\end{array}\right).
	\]

	Define
	\[
	\sqq :=\sqrt{q}.
	\]
	
	We continue to work with a hereditary abelian category $\ca$ as in \S\ref{subsec:periodic}.
	Let $\ch(\cc_1(\ca))$ be the Ringel-Hall algebra of $\cc_1(\ca)$ over $\Q(\sqq)$. 

	Following \cite{LP16,LW19a,LW20a}, we consider the ideal $\cI$ of $\ch(\cc_1(\ca))$ generated by
	\begin{align}
		\label{eq:ideal}
		&\Big\{ [M^\bullet]-[N^\bullet]\mid H(M^\bullet)\cong H(N^\bullet), \quad \widehat{\Im d_{M^\bullet}}=\widehat{\Im d_{N^\bullet}} \Big\}.
	\end{align}
	We denote
	\[
	\cs:=\{ a[K^\bullet] \in \ch(\cc_1(\ca))/\cI \mid a\in \Q(\sqq)^\times, K^\bullet\in \cc_1(\ca) \text{ acyclic}\},
	\]
	a multiplicatively closed subset of $\ch(\cc_1(\ca))/ \cI$ with the identity $[0]$.
	
	\begin{lemma}
		[\text{\cite[Proposition A.5]{LW19a}}]
		The multiplicatively closed subset $\cs$ is a right Ore, right reversible subset of $\ch(\cc_1(\ca))/\cI$. Equivalently, there exists the right localization of
		$\ch(\cc_1(\ca))/\cI$ with respect to $\cs$, denoted by $(\ch(\cc_1(\ca))/\cI)[\cs^{-1}]$.
	\end{lemma}

	The algebra $(\ch(\cc_1(\ca))/\cI)[\cs^{-1}]$ is the {\em semi-derived Ringel-Hall algebra} of $\cc_1(\ca)$ in the sense of \cite{LP16, LW19a} (also cf. \cite{Gor18}), and will be denoted by $\cs\cd\ch(\cc_1(\ca))$.

	\begin{definition}  \label{def:iH}
		The {\em $\imath$Hall algebra} of a hereditary abelian category $\ca$, denoted by $\iH(\ca)$,  is defined to be the twisted semi-derived Ringel-Hall algebra of $\cc_1(\ca)$, that is, the $\Q(\sqq)$-algebra on the same vector space as $\cs\cd\ch(\cc_1(\ca)) =(\ch(\cc_1(\ca))/\cI)[\cs^{-1}]$ equipped with the following modified multiplication (twisted via the restriction functor $\res: \cc_1(\ca)\rightarrow\ca$)
		\begin{align}
			\label{eq:tH}
			[M^\bullet]* [N^\bullet] =\sqq^{\langle \res(M^\bullet),\res(N^\bullet)\rangle_\ca} [M^\bullet]\diamond[N^\bullet].
		\end{align}
	\end{definition}

	For any $\alpha\in K_0(\ca)$,  there exist $X,Y\in\ca$ such that $\alpha=\widehat{X}-\widehat{Y}$. Define $[K_\alpha]:=[K_X]* [K_Y]^{-1}$. This is well defined, see e.g.,
	\cite[\S 3.2]{LP16}.
	It follows that $[K_\alpha]\; (\alpha\in K_0(\ca))$ are central in the algebra $\iH(\ca)$.
	
	
	The {\em quantum torus} $\widetilde{\ct}(\ca)$ is defined to be the subalgebra of $\iH(\ca)$ generated by $[K_\alpha]$, for $\alpha\in K_0(\ca)$. 
	
	\begin{proposition}
		[\text{\cite[Proposition 2.9]{LRW20a}}]
		\label{prop:hallbasis}
		The following hold in $\iH(\ca)$:
		\begin{enumerate}
			\item
			The quantum torus $\widetilde{\ct}(\ca)$ is a central subalgebra of $\iH(\ca)$.
			\item
			The algebra $\widetilde{\ct}(\ca)$ is isomorphic to the group algebra of the abelian group $K_0(\ca)$.
			\item
			$\iH(\ca)$ has an ($\imath$Hall) basis given by
			\begin{align*}
				\{[M]*[K_\alpha]\mid [M]\in\Iso(\ca), \alpha\in K_0(\ca)\}.
			\end{align*}
		\end{enumerate}
	\end{proposition}

	The following multiplication formula in the $\imath$Hall algebra is useful and applicable.
	\begin{proposition}
		[\text{\cite[Proposition 3.10]{LW20a}}]
		\label{prop:iHallmult}
		Let $\ca$ be a hereditary abelian category over $\bfk$.
		For any $A,B\in\ca\subset \cc_1(\ca)$, the following formulas hold in $\iH(\ca)$.
		\begin{align}
			\label{Hallmult1}
			[A]*[B]=&
			\sum_{[M]\in\Iso(\ca)}\sum_{[L],[N]\in\Iso(\ca)} \sqq^{-\langle A,B\rangle}  q^{\langle N,L\rangle}\frac{|\Ext^1(N, L)_{M}|}{|\Hom(N,L)| }
			\\
			\notag
			&\cdot |\{f\in\Hom(A,B)\mid \ker f\cong N, \coker f\cong L
			\}|[M]*[K_{\widehat{A}-\widehat{N}}],
			\\
			\label{Hallmult2}
			[A]*[B]=&
			\sum_{[M]\in\Iso(\ca)}\sum_{[L],[N]\in\Iso(\ca)} \sqq^{-\langle A,B\rangle}  \frac{|\Ext^1(N, L)_{M}|}{|\Ext^1(N,L)| }
			\\
			\notag
			&\cdot |\{f\in\Hom(A,B)\mid \ker f\cong N, \coker f\cong L
			\}|[M]*[K_{\widehat{A}-\widehat{N}}]
		\end{align}
	\end{proposition}

	\subsection{Realization of $\imath$quantum groups via quivers}
	\label{subsec:iquiver algebra}
	
	Let $Q=(Q_0,Q_1)$ be a quiver (not necessarily acyclic), and we sometimes write $\I =Q_0$.
	Let $\bfk=\F_{q}$ be a finite field of $q$ elements. 
	We denote by $\iH(\bfk Q)$ the twisted semi-derived Ringel-Hall algebra of $\cc_1(\rep_\bfk^{\rm nil}(Q))$. Denote by
	$\tUi_{|v= \sqq}$ or $\tUi_\sqq$ the $\imath$quantum group  specialised at $v=\sqq$.  
	Then we have the following result.

	\begin{theorem}
		[\text{\cite{LW20a}}]
		\label{thm:main}
		Let $Q$ be an arbitrary quiver. Then there exists a $\Q(\sqq)$-algebra embedding
		\begin{align*}
			\widetilde{\psi}_Q: \tUi_{|v= \sqq} &\longrightarrow \iH(\bfk Q),
		\end{align*}
		which sends
		\begin{align}
			B_i \mapsto \frac{-1}{q-1}[S_{i}],
			&\qquad
			\K_i \mapsto - q^{-1}[K_i], \qquad\text{ for }i \in \I.
			\label{eq:split}
		\end{align}
	\end{theorem}
	Let $\tCMH$  be the composition subalgebra of $\iH(\bfk Q)$ generated by $[S_i]$ and $[K_{S_i}]^{\pm1}$ for $i\in\I$.
	Then \eqref{eq:split} gives a $\Q(\sqq)$-algebra isomorphism
	\begin{align}
		\label{eq:isocomp}
		\widetilde{\psi}_Q: \tUi_{|v= \sqq} &\stackrel{\cong}{\longrightarrow} \tCMH.
	\end{align}

	For any sink $\ell \in Q_0$, define the quiver $\bs_\ell^+ Q$ by reversing all the arrows of $Q$ ending at $\ell $. 
	Associated to a sink $\ell \in Q_0$, the BGP reflection functor induces a {\rm reflection functor} (see \cite[\S3.2]{LW19b}):
	\begin{align}  \label{eq:Fl}
		F_\ell ^+: \cc_1(\rep_\bfk(Q))  \longrightarrow \cc_1(\rep_\bfk(\bs^+_\ell Q)).
	\end{align}
	The functor $F_\ell^+$ induces an isomorphism $\Gamma_\ell: \iH(\bfk Q)  \stackrel{\sim}{\rightarrow} \iH(\bfk \bs_\ell^+ Q)$ by \cite[Theorem 4.3]{LW19b},
	and then a commutative diagram
	\begin{equation}
		\label{eq:comm}
		\xymatrix{ \tUi_{ |v={\sqq}} \ar[r]^{ \TT_\ell} \ar[d]^{\widetilde{\psi}_{Q}} & \tUi_{ |v={\sqq}} \ar[d]^{\widetilde{\psi}_{\bs_\ell Q}}
			\\
			\iH(\bfk Q) \ar[r]^{\Gamma_\ell}  &\iH(\bfk \bs_\ell Q)}
	\end{equation}
	where $\TT_\ell$ is defined in Lemma \ref{lem:Ti}.	
	
	We consider the oriented cyclic quiver $C_n$ for $n\geq2$ with its vertex set $\Z_n=\{0,1,2,\dots,n-2,n-1\}$:
	\begin{center}\setlength{\unitlength}{0.5mm}
		\begin{equation}
			\label{fig:Cn}
			\begin{picture}(100,30)
				\put(2,0){\circle*{2}}
				\put(22,0){\circle*{2}}
				
				\put(82,0){\circle*{2}}
				\put(102,0){\circle*{2}}
				\put(52,25){\circle*{2}}
				
				\put(20,0){\vector(-1,0){16}}
				\put(40,0){\vector(-1,0){16}}
				\put(47.5,-2){$\cdots$}
				\put(80,0){\vector(-1,0){16}}
				\put(100,0){\vector(-1,0){16}}
				
				\put(54,24.5){\vector(2,-1){47}}
				\put(3,1){\vector(2,1){47}}
				
				\put(50.5,27){\tiny $0$}
				\put(1,-6){\tiny $1$}
				\put(21,-6){\tiny $2$}
				\put(75,-6){\tiny $n-2$}
				\put(95,-6){\tiny $n-1$}
			\end{picture}
		\end{equation}
		\vspace{-0.2cm}
	\end{center}
	Let $C_1$ be the Jordan quiver, i.e., the quiver with only one vertex and one loop arrow.
	
	Denote by $\rep^{\rm nil}_\bfk(C_n)$ the category of finite-dimensional nilpotent representations of $C_n$ over the field $\bfk$. Let $\iH(\bfk C_n)$ be the $\imath$Hall algebra of $\rep^{\rm nil}_\bfk(C_n)$.
	
	Define
	\begin{align}\label{the map Psi A}
		\Omega_{C_n}:=\widetilde{\psi}_{C_n}\circ \Phi:\tUiD_{|v=\sqq}\longrightarrow \iH(\bfk C_n).
	\end{align}
	Then $\Omega_{C_n}$ sends:
	\begin{align}
		\label{eq:HaDrA1}\K_j\mapsto [K_{S_j}],\qquad C\mapsto [K_\de],\qquad
		B_{j,0}\mapsto \frac{-1}{q-1}[S_j], \qquad \forall \,\,1\leq j\leq n-1.
	\end{align}
	For any $1\leq j\leq n-1$, $l\in\Z$ and $r\geq 1$, we define
	\begin{align}
		\label{def:haB}
		\haB_{j,l}:=(1-q)\Omega_{C_n}(B_{j,l}),\qquad
		\widehat{\Theta}_{j,r}:= \Omega_{C_n}(\Theta_{j,r}),\qquad
		\widehat{H}_{j,r}:= \Omega_{C_n}(H_{j,r}).
	\end{align}
	In particular, $\haB_{j,0}=[S_j]$ for any $1\leq j\leq n-1 $.
	
	It is helpful and interesting to compute all root vectors 
	clearly in $\iH(\bfk C_n)$.  There is a way to interpret the root vectors in $\imath$Hall algebras, which we give in the following.
	
	Let $Q$ be a quiver with its underlying graph the same as $C_n$. Then there exists an isomorphism $\texttt{FT}_{C_n,Q}:\tCMH\stackrel{\cong}{\longrightarrow} \tCMHC$ 
	such that $\texttt{FT}_{C_n,Q}([S_i])=[S_i]$, $\texttt{FT}_{C_n,Q}([K_{S_i}])=[K_{S_i}]$ for any $i\in\I$; cf. \eqref{eq:isocomp}. For any positive real root $\beta$ of $\widehat{\mathfrak{sl}}_n$, we denote by $M(\beta)$ the unique (up to isomorphism) $\bfk Q$-module with its class $\beta$ in $K_0(\mod(\bfk Q))$.

	For any hereditary abelian category $\ca$ and a full subcategory $\cs$, recall from \cite{GL87} that the right (resp. left) perpendicular category $\cs{}^\bot$ (resp. ${}^\bot \cs$) is defined as follows:
	\begin{align}
		\label{def:perpcat right}
		\cs{}^\bot:=\{X\in\ca\mid \Hom_\ca(C,X)=0=\Ext_\ca^1(C,X)\text{ for any }C\in\cs\};\\
		\label{def:perpcat}
		{}^\bot\cs:=\{X\in\ca\mid \Hom_\ca(X,C)=0=\Ext_\ca^1(X,C)\text{ for any }C\in\cs\}.
	\end{align}
	Then the category $\cs{}^\bot$ (resp. ${}^\bot \cs$) is a hereditary abelian category, which is extension-closed in $\ca$; see \cite[Proposition 3.2]{BS13}.
	If $\cs=\{M\}$ for an object $M\in\ca$, then we also denote $\cs{}^\bot:=M{}^\bot$ and ${}^\bot\cs:={}^\bot M$.

	Fix $1\leq j\leq n-1$. We assume $Q(j)$ is the following quiver
	\begin{center}\setlength{\unitlength}{0.5mm}
		\begin{equation}
			\label{quiver:affineA}
			\begin{picture}(100,30)
				\put(2,0){\circle*{2}}
				\put(51,0){\circle*{2}}
				
				\put(102,0){\circle*{2}}
				\put(52,25){\circle*{2}}
				
				\put(3.5,0){\vector(1,0){16}}
				\put(21,-2){$\cdots$}
				
				\put(32,0){\vector(1,0){16}}
				
				\put(71.5,-2){$\cdots$}
				\put(70,0){\vector(-1,0){16}}
				\put(100,0){\vector(-1,0){16}}
				
				\put(54,24.5){\vector(2,-1){47}}
				\put(50,24.5){\vector(-2,-1){47}}
				
				\put(50.5,27){\tiny $0$}
				\put(1,-6){\tiny $1$}
				\put(50,-6){\tiny $j$}
				\put(95,-6){\tiny $n-1$}
			\end{picture}
		\end{equation}
		\vspace{-0.2cm}
	\end{center}

	Recall that
	\begin{align}
		\label{eq:wj}
		\omega_j=\sigma^j (s_{n-j}\cdots s_1)(s_{n-j+1}\cdots s_2) \cdots (s_{n-1}\cdots s_j),
	\end{align}
	where $\sigma$ is the automorphism of the underlying graph of $C_n$ such that $\sigma(i)=i+1$; see e.g. \cite[\S6.2]{DJX12}.
	Note that $\bs_{\omega_j}(Q(j))=Q(j)$, and $j,\dots, n-1,\dots, 2,\dots, n-j+1,1,\dots, n-j$  is an $\imath$-admissible sequence. 
	
	\begin{proposition}
		\label{rem:rootvectors}
		For any $l>0$, we have
		\begin{align}
			\label{eq:haB-}
			\haB_{j,-l}=& (-1)^{jl}\texttt{FT}_{C_n,Q(j)}([M(l\de-\alpha_j)])*[K_{l\de-\alpha_j}]^{-1},
			\\
			\label{eq:haB+}
			\haB_{j,l}= &(-1)^{jl}\texttt{FT}_{C_n,Q(j)}([M(l\de+\alpha_j)]),
			\\
			\label{eq:haTh+}
			\widehat{\Theta}_{j,l}=&\frac{(-1)^{jl}}{(q-1)^2\sqq^{l-1}}\texttt{FT}_{C_n,Q(j)}\Big(\sum_{0\neq f:S_j\rightarrow M(l\de+\alpha_j) } [\coker f]\Big).
		\end{align}
	\end{proposition}
	
	\begin{proof}
		From the proof of \cite[Proposition 5
		.1]{LR21}, in $\iH(\bfk Q(j))$ we have
		\begin{align*}
			\Gamma_{\omega_j}^l([S_j])=[M(l\de-\alpha_j)]*[K_{l\de-\alpha_j}]^{-1}.
		\end{align*}
		Then $B_{j,-l}$ satisfies that $\Phi(B_{j,-l})= (-1)^{jl} \TT_{\omega_j}^{l} (B_j)$ and so
		$$\widetilde{\psi}_{Q(j)}\circ\Phi(B_{j,-l})=(-1)^{jl}\cdot\frac{-1}{q-1}[M(l\de-\alpha_j)]*[K_{l\de-\alpha_j}]^{-1} $$
		in $\iH(\bfk Q(j))$. Similarly, one can see that
		$$\widetilde{\psi}_{Q(j)}\circ\Phi(B_{j,l})=(-1)^{jl}\cdot\frac{-1}{q-1}[M(l\de+\alpha_j)];$$
		cf. \cite[Lemma 5.5]{LRW20a}.
		So we have \eqref{eq:haB-}--\eqref{eq:haB+}
		
		For \eqref{eq:haTh+}, let $\cs:=\{S_i\in\mod(\bfk Q(j))\mid i\neq 0,j\}$. Then the subcategory ${}^\bot\cs$ is equivalent to $\mod(\bfk Q_{\texttt{Kr}})$, where $Q_{\texttt{Kr}}$ is the Kronecker quiver $\xymatrix{0\ar@<0.5ex>[r] \ar@<-0.5ex>[r] & 1}$. So there exists an embedding functor $\F: \mod(\bfk Q_{\texttt{Kr}})\simeq {}^\bot\cs\rightarrow \mod(\bfk Q(j))$. In particular, $\F(S_1)=S_j$, and $\F(M(l\alpha_0+(l+1)\alpha_1))=M(l\de+\alpha_j)$ for any $l\geq0$.
		So $\F$ induces an embedding of algebras $\widetilde{\ch}(\cc_1(\mod(\bfk Q_{\texttt{Kr}})))\rightarrow \widetilde{\ch}(\cc_1(\mod(\bfk Q(j))))$, and then an embedding of algebras $F:\iH(\bfk Q_{\texttt{Kr}})\rightarrow \iH(\bfk Q (j))$.
		Using \cite[Theorem 5.11]{LRW20a}, in order to compute $\haTh_{j,l}$ in $\iH(\bfk Q_{\texttt{Kr}})$, it is equivalent to compute them in $\iH(\P^1_\bfk)$. By  \cite[(4.2)]{LRW20a}, we have
		$$\haTh_{1,l}=\frac{1}{(q-1)^2\sqq^{l-1}}\sum_{0\neq f:S_1\rightarrow M(l\alpha_0+(l+1)\alpha_1) } [\coker f]$$
		in $\iH(\bfk Q_{\texttt{Kr}})$. By applying $F$, we have
		\begin{align}
			\label{haTh}
			\haTh_{j,l}=\frac{(-1)^{jl}}{(q-1)^2\sqq^{l-1}}\sum_{0\neq f:S_j\rightarrow M(l\de+\alpha_j) } [\coker f]
		\end{align}
		in $\iH(\bfk Q (j))$. Here the scalar $(-1)^{jl}$ comes from \eqref{eq:haB-}--\eqref{eq:haB+}. Then \eqref{eq:haTh+} follows.
	\end{proof}

	\subsection{Weighted projective lines}
	Recall $\bfk=\F_q$. Fix a positive integer $\bt$ such that $2\leq \bt\leq q$ in the following. Let $\X=\X_\bfk=\X_{\bp,\ul{\bla}}$ be the weighted projective line attached to 
	a weight sequence $\bp=(p_1,p_2,\dots,p_\bt)\in\Z_{+}^{\bt}$, and a parameter sequence $\ul{\bla}=\{\bla_1,\dots,\bla_{\mathbf{t}}\}$ of distinguished closed points (of degree one) on the projective line $\PL$. Let $L(\bp)$
	be the rank one abelian group on generators $\vec{x}_1$, $\vec{x}_2$, $\cdots$, $\vec{x}_{\bt}$ with relations
	$p_1\vec{x}_1=p_2\vec{x}_2=\cdots=p_{\bt}\vec{x}_{\bt}$. We call $\vec{c}:=p_i\vec{x}_i$ the canonical element of $L(\bp)$.
	
	Let $\co:=\co_\X$ be the  structure sheaf, and $\coh(\X)$ the category of coherent sheaves \cite{GL87}.

	Let $\scrf$ be the full subcategory of $\coh(\X)$ consisting of all locally free sheaves, and $\scrt$ be the full subcategory consisting of all torsion sheaves. 

	Then the structure of $\scrt$ is described via the representation theory of cyclic quivers.
	\begin{lemma}[\cite{GL87}; also see \cite{Sch04}]
		\label{lem:isoclasses Tor}
		(1) The category $\scrt$ decomposes as a coproduct $\scrt=\coprod_{x\in\X} \scrt_{x}$, where $\scrt_{x}$ is the subcategory of torsion sheaves with support at $x$.
		
		(2) For any ordinary point $x$ of degree $d$, let $\bfk_{x}$ denote the residue field at $x$, i.e., $[\bfk_x:\bfk]=d$. Then $\scrt_{x}$ is equivalent to the category $\rep^{\rm nil}_{\bfk_{x}}(C_1)$.
		
		(3) For any exceptional point $\bla_i$ ($1\leq i\leq \bt$), the category $\scrt_{\bla_i}$ is equivalent to $\rep^{\rm nil}_{\bfk}(C_{p_i})$.
	\end{lemma}

	For any ordinary point $\blx$ of degree $d$, let $\pi_{\blx}$ be the prime homogeneous polynomial corresponding to $\blx$. Then multiplication by $\pi_{\blx}$ gives the exact sequence
	$$0\longrightarrow \co\stackrel{\pi_{\blx}}{\longrightarrow} \co(d\vec{c})\longrightarrow S_{\blx}\rightarrow0,$$
	where $S_{\blx}$ is the unique (up to isomorphism) simple sheaf in the category $\scrt_{\blx}$. Then $S_{\blx}(\vec{l})=S_{\blx}$ for any $\vec{l}\in L(\bp)$.
	
	For any exceptional point $\bla_i$, multiplication by $x_i$ yields the short exact sequence
	$$0\longrightarrow \co((j-1)\vec{x}_i)\stackrel{x_i}{\longrightarrow} \co(j\vec{x}_i)\longrightarrow S_{ij}\rightarrow0,\text{ for } 1\le j\le p_i;$$
	where $\{S_{ij}\mid  j\in\Z_{p_i}\}$ is a complete set of pairwise non-isomorphic simple sheaves in the category $\scrt_{\bla_i}$ for any $1\leq i\leq \bt$. 

	We denote by $S_{\blx}^{(a)}$ the unique indecomposable object of length $a$ in $\scrt_{\blx}$ for any ordinary point ${\blx}$; and denote by $S_{ij}^{(a)}$ the unique indecomposable object with top $S_{ij}$ and length $a$ in $\scrt_{\bla_i}$.

	The Grothendieck group $K_0(\coh(\X))$ of $\coh(\X)$ satisfies
	\begin{align}
		K_0(\coh(\X))\cong \Big(\Z\widehat{\co}\oplus \Z\widehat{\co(\vec{c})}\oplus \bigoplus_{i,j}\Z\widehat{S_{ij}}\Big)/I,
	\end{align}
	where $I$ is the subgroup generated by $\{\sum_{j=1}^{p_i} \widehat{S_{ij}} +\widehat{\co} -\widehat{\co(\vec{c})}\mid i=1,\dots,\bt\}$; see \cite{GL87}.
	Let $\de:=\widehat{\co(\vec{c})}-\widehat{\co}= \sum_{j=1}^{p_i} \widehat{S_{ij}}$.
	
	
	Recall the root system $\cR$ of $L(\fg)$. Then there is a natural isomorphism of $\Z$-modules $K_0(\coh(\X))\cong \cR$ given as below (for $1\leq i\leq \bt,\;1\le j\le p_i-1, $ $r\in\N$ and $l\in\Z$):
	\begin{align*}
		&\widehat{S_{ij}} \mapsto \alpha_{ij}, 
		&\widehat{S_{i,0}}\mapsto \de-\sum_{j=1}^{p_i-1}\alpha_{ij},\quad
		&\widehat{S_{i,0}^{(p_i-1)}}\mapsto \de-\alpha_{i1}, 
		&\widehat{S_{i,0}^{(rp_i)}}\mapsto r\de,\quad
		&\widehat{\co(l\vec{c})}\mapsto \alpha_\star+l\de.
	\end{align*}
	Under this isomorphism, the symmetric Euler form on $K_0(\coh(\X))$ coincides with the Cartan form on $\cR$. So we always identify $K_0(\coh(\X))$ with $\cR$ via this isomorphism in the following.

	\section{$\imath$Hall algebras and $\imath$quantum loop algebras}
	\label{sec:hom}
	
	In this section, we shall formulate the main result of this paper.
	
	\subsection{Root vectors}
	\label{subsec:embeddingtube}
	Let $\X$ be a weighted projective line of weight type $(\bp,\ul{\bla})$.
	Recall that $\bla_i$ is the exceptional closed point of $\X$ of weight $p_i$ for any $1\leq i\leq \bt$, and $\scrt_{\bla_i}$ is the Serre subcategory of $\coh(\X)$ consisting of torsion sheaves supported at $\bla_i$. There is an  embedding of $\imath$Hall algebras:
	\begin{equation}
		\label{eq:embeddingx}
		\iota_i: \iH(\bfk C_{p_i})\longrightarrow\iH(\X_\bfk).
	\end{equation}
	
	Inspired by \eqref{def:haB}, we define
	\begin{align}
		\label{def:haBThH}
		\haB_{[i,j],l}:= \iota_i (\haB_{j,l}),\quad \widehat{\Theta}_{[i,j],r}:=\iota_i(\widehat{\Theta}_{j,r}),\quad \widehat{H}_{[i,j],r}:=\iota_{i}(\widehat{H}_{j,r}),
	\end{align}
	for any $1\leq j\leq p_i-1$, $l\in\Z$ and $r>0$.
	By definition and \eqref{exp h}, we have
	\begin{align}
		\label{exp haH}
		1+ \sum_{m\geq 1} (\sqq-\sqq^{-1})\haTh_{[i,j],m} u^m  = \exp\Big( (\sqq-\sqq^{-1}) \sum_{m\geq 1} \widehat{H}_{[i,j],m} u^m \Big).
	\end{align}
	
	For convenience, we set $\haTh_{[i,j],0}= \frac{1}{\sqq-\sqq^{-1}}$ for any $1\leq i\leq \bt,1\le j\le p_i-1$.

	
	
	
	
	Let $\mathcal{C}$ be the Serre subcategory of $\coh(\X)$ generated by those simple sheaves $S$ satisfying $\Hom(\co, S)=0$. Recall from \cite{BS13} that the Serre quotient $\coh(\X)/\mathcal{C}$ is equivalent to the category $\coh(\P^1)$ and the canonical functor $\coh(\X)\rightarrow\coh(\X)/\mathcal{C}$ has an exact fully faithful right adjoint functor
	\begin{equation}
		\label{the embedding functor F}
		\mathbb{F}_{\X,\P^1}: \coh(\P^1)\rightarrow\coh(\X),
	\end{equation}
	which sends $$\co_{\P^1}(l)\mapsto\co(l\vec{c}),\quad S_{\bla_i}^{(r)}\mapsto S_{i,0}^{(rp_i)},\quad S_{x}^{(r)}\mapsto S_{x}^{(r)}$$ 
	for any $l\in\Z, r\geq 1$, $1\leq i\leq \bt$ and $x\in\PL\setminus\{\bla_1,\cdots, \bla_\bt\}$. Then $\mathbb{F}_{\X,\P^1}$ induces an exact fully faithful functor $\cc_1(\coh(\P^1))\rightarrow \cc_1(\coh(\X))$, which is also denoted by $\mathbb{F}_{\X,\P^1}$. This functor $\mathbb{F}_{\X,\P^1}$ induces a canonical embedding of Hall algebras $\ch\big(\cc_1(\coh(\P^1))\big)\rightarrow \ch\big(\cc_1(\coh(\X))\big)$, and then an embedding of $\imath$Hall algebras:
	\begin{equation}
		\label{the embedding functor F on algebra}
		F_{\X,\P^1}:\iH(\P^1_\bfk) \longrightarrow \iH(\X_\bfk).
	\end{equation}
	
	Recall from \cite[(4.2), Proposition 6.3]{LRW20a} that for any $m\geq0$, in $\iH(\P^1_\bfk)$
	$$\haT_{m}= \frac{1}{(q-1)^2\sqq^{m-1}}\sum_{0\neq f:\co(s)\rightarrow \co(m+s) } [\coker f],
	$$
	which is independent of $s\in\Z$, and
	\begin{align*}
		\haH_{m}
		=&\sum_{x,d_x|m} \frac{[m]_\sqq}{m} d_x \sum_{|\lambda|=\frac{m}{{d_x}}} \bn_x(\ell(\lambda)-1)\frac{[S_x^{(\lambda)}]}{\big|\Aut(S_x^{(\lambda)})\big|}  -\de_{m,ev}\frac{[m]_\sqq}{m} [K_{\frac{m}{2}\de}],
	\end{align*}
	where $d_x$ denotes the degree of $x\in\P^1_\bfk$, and
	\[\bn_x(l)=\prod_{i=1}^l(1-\sqq^{2id_x}), \,\,\forall l\geq1.\]

	Hence, we define
	\begin{align}
		\label{def:Theta star}
		\widehat{\Theta}_{\star,m}:=F_{\X,\P^1}(\haT_m)=\frac{1}{(q-1)^2\sqq^{m-1}}\sum_{0\neq f:\co(s\vec{c})\rightarrow \co((m+s)\vec{c}) } [\coker f],
	\end{align}
	and
	\begin{align}
		\label{formula for Hm}
		\haH_{\star,m}:=F_{\X,\P^1}(\haH_{m})
		=&\sum_{x,d_x|m} \frac{[m]_\sqq}{m} d_x \sum_{|\lambda|=\frac{m}{{d_x}}} \bn_x(\ell(\lambda)-1)\frac{F_{\X,\P^1}([S_x^{(\lambda)}])}{\big|\Aut(S_x^{(\lambda)})\big|}  -\de_{m,ev}\frac{[m]_\sqq}{m} [K_{\frac{m}{2}\de}].
	\end{align}
	Here, for any partition $\lambda=(\lambda_1, \lambda_2,\cdots)$, $F_{\X,\P^1}([S_x^{(\lambda)}])=[S_x^{(\lambda)}]$ for $x\in\PL\setminus\{\bla_1,\cdots, \bla_\bt\}$, while $F_{\X,\P^1}([S_{\bla_i}^{(\lambda)}])=S_{i,0}^{(\lambda)}:=\bigoplus_{i}S_{i,0}^{(\lambda_ip_i)}$ for $1\leq i\leq \bt$.
	Note that $\haTh_{\star,m}=0$ if $m<0$ and $\haTh_{\star,0}=\frac{1}{\sqq-\sqq^{-1}}$.

	\begin{lemma}[\text{\cite[Lemma 6.1]{LR21}}]
		\label{lem:Hxm}	
		For any $m\geq 1$ and $x\in\P^1_\bfk$ such that $d_x| m$, we define
		\begin{align}
			\label{def:Hxm}
			\widehat{H}_{x,m}:= \frac{[m]_\sqq}{m} d_x \sum_{|\lambda|=\frac{m}{{d_x}}} \bn_x(\ell(\lambda)-1)\frac{F_{\X,\P^1}([S_x^{(\lambda)}])}{\big|\Aut(S_x^{(\lambda)})\big|}  -\de_{\frac{m}{d_x},ev}\cdot d_x\sqq^{-\frac{m}{2}}\frac{[\frac{m}{2}]_\sqq }{m} [K_{\frac{m}{2}\de}].
		\end{align}
		Then
		\begin{align}
			\label{eq:Hdecomptox}
			\widehat{H}_{\star,m}=\sum_{x,d_x|m} \widehat{H}_{x,m}.
		\end{align}
	\end{lemma}

	
	
	\subsection{The homomorphism $\Omega$}
	\label{subsec:homo}
	
	Recall the star-shaped graph $\Gamma=T_{p_1,\dots,p_\bt}$ in \eqref{star-shaped}. The following is one of the main results of this paper.
	
	\begin{theorem}
		\label{thm:morphi}
		For a star-shaped graph $\Gamma$ of finite or affine type, let $\fg$ be the Kac-Moody algebra and $\X$ be the weighted projective line associated to $\Gamma$. Then there exists a $\Q(\sqq)$-algebra embedding
		\begin{align}
			\Omega: \tUiD_{ |v=\sqq}\longrightarrow \iH(\X_\bfk),
		\end{align}
		which sends
		\begin{align}
			\label{eq:mor1}
			&\K_{\star}\mapsto [K_{\co}], \qquad \K_{[i,j]}\mapsto [K_{S_{ij}}], \qquad C\mapsto [K_\de];&
			\\
			\label{eq:mor2}
			&{B_{\star,l}\mapsto \frac{-1}{q-1}[\co(l\vec{c})]},\qquad\Theta_{\star,r} \mapsto {\widehat{\Theta}_{\star,r}}, \qquad H_{\star,r} \mapsto {\widehat{H}_{\star,r}};\\
			\label{eq:mor3}
			&\y_{[i,j],l}\mapsto {\frac{-1}{q-1}}\haB_{[i,j],l}, \quad \Theta_{[i,j],r}\mapsto\widehat{\Theta}_{[i,j],r}, \quad H_{[i,j],r}\mapsto \widehat{H}_{[i,j],r};
		\end{align}
		for any $[i,j]\in\II-\{\star\}$, $l\in\Z$, $r>0$.
	\end{theorem}
	
	By \cite{LR21}, we already know the morphism $\Omega: \tUiD_{ |v=\sqq}\rightarrow \iH(\X_\bfk)$ is  well defined, and we only need to prove that $\Omega$ is injective, which shall occupy the following two sections.

	The composition subalgebra $\tCMHX$ is the subalgebra of $\tMHX$ generated by the elements $[\co(l\vec{c})]$ for $l\in\Z$, $\widehat{\Theta}_{\star,r}$ for $r\geq1$, $[S_{ij}]$ for $1\leq i\leq \bt$ and $1\leq j\leq p_i-1$, together with $[K_{\co}]^{\pm1}, [K_{\delta}]^{\pm1}$ and $[K_{S_{ij}}]^{\pm1}$. Then the following corollary follows from Theorem \ref{thm:morphi} and Lemma \ref{lem:reduced generators} immediately.
	
	\begin{corollary}
		\label{cor:epimorphism}
		For a star-shaped graph $\Gamma$ of finite or affine type, let $\fg$ be the Kac-Moody algebra and $\X$ be the weighted projective line associated to $\Gamma$. Then there exists a $\Q(\sqq)$-algebra isomorphism
		\begin{align}
			\Omega: \tUiD_{ |v=\sqq}\longrightarrow \tCMHX,
		\end{align}
		defined as in \eqref{eq:mor1}--\eqref{eq:mor3}.
	\end{corollary}

	\begin{remark}
		The structure of $\tCMHX$ is rich, which will be studied in depth in a separate publication; compare with \cite{BS13}. Moreover, a PBW basis will be given for the $\imath$quantum group via coherent sheaves if $\X_\bfk$ is of domestic type.
	\end{remark}
	
	
	%



	\section{$\imath$Hall algebras of cyclic quivers and $\tUi_v(\widehat{\mathfrak{gl}}_n)$}
	\label{sec:Hall cyclic}
	
	In this section, we shall investigate the $\imath$Hall algebra of a cyclic quiver $C_n$ again, this is a preparation of proving the injectivity of $\Omega:\tUiD_{|v=\sqq}\rightarrow \iH(\X_\bfk)$ obtained in Theorem \ref{thm:morphi} for $\fg$ of finite or affine type. As a byproduct, we use the whole $\iH(\bfk C_n)$ to realize the $\imath$quantum loop algebra of $\widehat{\mathfrak{gl}}_n$.

	\subsection{Ringel-Hall algebras and quantum groups}
	
	Let $Q$ be a quiver (not necessarily acyclic) with vertex set $Q_0= \I$.
	Let $n_{ij}$ be the number of edges connecting vertex $i$ and $j$. Let $C=(c_{ij})_{i,j \in \I}$ be the symmetric generalized Cartan matrix of the underlying graph of $Q$, defined by $c_{ij}=2\delta_{ij}-n_{ij}.$ Let $\fg$ be the corresponding Kac-Moody algebra.
	The Drinfeld double quantum group $\tU=\tU_v(\fg)$ \cite{LW19a} is defined to be the $\Q(\bv)$-algebra generated by $E_i,F_i, \tK_i,\tK_i'$, $i\in \I$, where $\tK_i, \tK_i'$ are invertible, subject to the following relations:
	\begin{align}
		[E_i,F_j]= \delta_{ij} \frac{\tK_i-\tK_i'}{\bv-\bv^{-1}},  &\qquad [\tK_i,\tK_j]=[\tK_i,\tK_j']  =[\tK_i',\tK_j']=0,
		\label{eq:KK}
		\\
		\tK_i E_j=\bv^{c_{ij}} E_j \tK_i, & \qquad \tK_i F_j=\bv^{-c_{ij}} F_j \tK_i,
		\label{eq:EK}
		\\
		\tK_i' E_j=\bv^{-c_{ij}} E_j \tK_i', & \qquad \tK_i' F_j=\bv^{c_{ij}} F_j \tK_i',
		\label{eq:K2}
	\end{align}
	and the quantum Serre relations, for $i\neq j \in \I$,
	\begin{align}
		& \sum_{r=0}^{1-c_{ij}} (-1)^r \left[ \begin{array}{c} 1-c_{ij} \\r \end{array} \right]  E_i^r E_j  E_i^{1-c_{ij}-r}=0,
		\label{eq:serre1} \\
		& \sum_{r=0}^{1-c_{ij}} (-1)^r \left[ \begin{array}{c} 1-c_{ij} \\r \end{array} \right]  F_i^r F_j  F_i^{1-c_{ij}-r}=0.
		\label{eq:serre2}
	\end{align}

	Let $\widetilde{\bU}^+$ be the subalgebra of $\widetilde{\bU}$ generated by $E_i$ $(i\in \I)$, $\widetilde{\bU}^0$ be the subalgebra of $\widetilde{\bU}$ generated by $\tK_i, \tK_i'$ $(i\in \I)$, and $\widetilde{\bU}^-$ be the subalgebra of $\widetilde{\bU}$ generated by $F_i$ $(i\in \I)$, respectively.

	Let $\ch(\bfk Q)$ be the  Ringel-Hall algebra of $\rep_\bfk^{\rm nil}(Q)$; see \eqref{eq:mult}.
	Let $\widetilde{\ch}(\bfk Q)$ be the twisted Ringel-Hall algebra of $\rep_\bfk^{\rm nil}(Q)$, that is,
	\begin{align}
		\label{eq:RHallmult}
		[M]* [N]=\sqq^{\langle M,N\rangle_{Q}}[M]\diamond [N], \quad \forall M,N\in \rep_\bfk^{\rm nil}(Q).
	\end{align}
	From \cite{Rin90,Gr95}, there exists an algebra embedding
	\begin{align}
		\label{eq:Xi}
		\Xi_{Q}:\tU^+|_{v=\sqq}&\longrightarrow \widetilde{\ch}(\bfk Q)
		\\
		E_i&\mapsto\frac{-1}{q-1}[S_i], \text{ for any }i\in\I.
	\end{align}
	
	Let $\widetilde{\cc}(\bfk Q)$ be the composition subalgebra of $\widetilde{\ch}(\bfk Q)$ generated by $[S_{i}]$ ($i\in\I$). Then $\Xi_Q$ induces an algebra isomorphism $\tU^+|_{v=\sqq}\stackrel{\cong}{\longrightarrow} \widetilde{\cc}(\bfk Q)$.

	\subsection{Semi-derived Ringel-Hall algebras of double framed quiver algebras}
	Let $Q=(Q_0,Q_1)$ be a quiver. Let $R_2$ be the radical square zero path algebra of $\xymatrix{1' \ar@<0.5ex>[r]^{\varepsilon'} & 1 \ar@<0.5ex>[l]^{\varepsilon}}$, i.e., $\varepsilon' \varepsilon =0 =\varepsilon\varepsilon '$.
	Consider $\Lambda=\bfk Q\otimes_\bfk R_2$ defined as in \cite{LW19a}, which is called the double framed quiver algebra of $Q$.
	By \cite[Section 2]{LW19a}, we have
	\begin{align}
		\Lambda\cong \bfk Q^\sharp/I^\sharp,
	\end{align}
	where
	$Q^{\sharp}$ is the quiver such that
	\begin{itemize}
		\item[-] the vertex set of $Q^{\sharp}$ consists of 2 copies of the vertex set $Q_0$, $\{i,i'\mid i\in Q_0\}$;
		\item[-] the arrow set of $Q^{\sharp}$ is
		\[
		\{\alpha: i\rightarrow j,\;\alpha': i'\rightarrow j'\mid\alpha\in Q_1\}\cup\{ \varepsilon_i: i\rightarrow i' ,\;\varepsilon'_i: i'\rightarrow i\mid i\in Q_0 \};
		\]
	\end{itemize}
	and $I^\sharp$
	is the two-sided ideal of $\bfk Q^{\sharp}$ generated by
	\begin{itemize}
		\item[-]
		(Nilpotent relations) $\varepsilon_i \varepsilon_i'$, $\varepsilon_i'\varepsilon_i$ for any $i\in Q_0$;
		\item[-]
		(Commutative relations) $\varepsilon_j' \alpha' -\alpha\varepsilon_i'$, $\varepsilon_j \alpha -\alpha'\varepsilon_i$ for any $(\alpha:i\rightarrow j)\in Q_1$.
	\end{itemize}
	
	Let $Q'$ be an identical copy of $Q$. Then $Q^{\dbl}:=Q\sqcup Q'$ is a subquiver of $Q^\sharp$. Note that $Q^{\dbl}$ admits a natural involution $\swa$. As a consequence, $\Lambda$ is the $\imath$quiver algebra of the $\imath$quiver $(Q^{\dbl},\swa)$ in the sense of \cite{LW19a}; see \cite[Example 2.3]{LW19a}.
	
	Let $\mod^{\rm nil}(\Lambda)$ be the category of finite-dimensional nilpotent $\Lambda$-modules. Then $\mod^{\rm nil}(\Lambda)$ is equivalent to the category $ \cc_{\Z/2}(\rep_{\bfk}^{\rm nil}(Q))$ of $2$-periodic complexes over $\rep_{\bfk}^{\rm nil}(Q)$; see \cite{Br13,LP16}.
	Let $\cp^{<\infty}(\Lambda)$ be the subcategory of $\mod^{\rm nil}(\Lambda)$ formed by modules with finite projective dimension.
	
	For any $i\in Q_0$, let
	\begin{align}
		\E_i:=\xymatrix{\bfk\ar@<0.5ex>[r]^1 & \bfk\ar@<0.5ex>[l]^0} \text{ and } \E_{i'}:=\xymatrix{\bfk\ar@<0.5ex>[r]^0 & \bfk\ar@<0.5ex>[l]^1} \text{ on the quiver } \xymatrix{i'\ar@<0.5ex>[r]^{\varepsilon_{i'}} & i\ar@<0.5ex>[l]^{\varepsilon_{ i}} }.
	\end{align}
	Then $\E_i,\E_{i'}$ can be viewed as $\Lambda$-modules. By \cite[Lemma 2.3]{LW20a}, we have $\E_i,\E_{i'}\in \cp^{<\infty}(\Lambda)$.
	
	Let $\widetilde{\ch}(\Lambda)$ be the twisted Hall algebra of $\mod^{\rm nil}(\Lambda)$ with the product
	\begin{align}
		\label{eqn:twsited multiplication}
		[M]* [N] =\sqq^{\langle \res(M),\res(N)\rangle_{Q\sqcup Q'}} [M]\diamond[N]
	\end{align}
	for any $M,N\in \mod^{\rm nil}(\Lambda)$. Here $\res:\mod^{\rm nil}(\Lambda)\rightarrow \rep_\bfk^{\rm nil}(Q\sqcup Q')$ is the restriction functor by viewing $\bfk (Q\sqcup Q')$ as a subalgebra of $\Lambda$; see \cite[\S2]{LW19a}.
	
	Let $I$ be the two-sided ideal of $\widetilde{\ch}(\Lambda)$ generated by all the differences $[L]-[K\oplus M]$ if there is a short exact sequence
	\begin{equation}
		\label{eq:ideal}
		0 \longrightarrow K \longrightarrow L \longrightarrow M \longrightarrow 0
	\end{equation}
	in $\mod^{\rm nil}(\Lambda)$  with $K\in \cp^{<\infty}(\Lambda)$.
	Let $\cs_\Lambda$ be the multiplicatively closed subset of the quotient algebra $\widetilde{\ch}(\Lambda)/I $:
	\begin{equation}
		\label{eq:Sca}
		\cs_{\Lambda} := \{ a[K] \in \widetilde{\ch}(\Lambda)/I \mid a\in \Q(\sqq)^\times, K\in \cp^{<\infty}(\Lambda)\}.
	\end{equation}
	
	Note that $\Lambda$ is $1$-Gorenstein.
	Following \cite{LP16,LW19a}, we define the twisted semi-derived Ringel-Hall algebra of $\Lambda$ to be
	the right localization of
	$\widetilde{\ch}(\Lambda)/I$ with respect to $\cs_{\Lambda}$, and denoted by $\cs\cd\widetilde{\ch}(\Lambda)$ or $\iH(\bfk Q^{\dbl},\swa)$ in the following.
	
	Note that $\bfk Q$ and $\bfk Q'$ are quotient algebras of $\Lambda$. We view the twisted Ringel-Hall algebra $\widetilde{\ch}(\bfk Q)$ and $\widetilde{\ch}(\bfk Q')$ as subalgebras of $\cs\cd\widetilde{\ch}(\Lambda)$; see \cite[\S8]{LW19a}.
	
	From \cite[Theorem 4.9]{LP16}, we see that $\cs\cd\widetilde{\ch}(\Lambda)$ is isomorphic to the Drinfeld double of the twisted extended Ringel-Hall algebra of $\rep_{\bfk}^{\rm nil}(Q)$. The following result is a generalization of Bridgeland's result (without assuming $Q$ to be acyclic).

	\begin{lemma}[\text{cf. \cite[Theorem 8.5]{LW19a}}]
		\label{lem:UslnSDH}
		Let $Q$ be a quiver. Then there exists an algebra embedding
		\begin{align}
			\label{eq:semiXi}
			\Xi_Q:\tU|_{v=\sqq}&\longrightarrow \cs\cd\widetilde{\ch}(\Lambda)
			\\
			E_i\mapsto&\frac{-1}{q-1}[S_i], \qquad F_i\mapsto  \frac{1}{\sqq-\sqq^{-1}}[S_{i'}]
			\\
			\tK_i\mapsto&[\E_{i}], \qquad \qquad \tK_i'\mapsto[\E_{i'}],
			\qquad\forall i\in\I.
		\end{align}
	\end{lemma}
	
	\begin{proof}
		Since $\cs\cd\widetilde{\ch}(\Lambda)$ is isomorphic to the Drinfeld double of the twisted extended Hall algebra of $\rep_\bfk^{\rm nil}(Q)$, the result follows from
		\eqref{eq:Xi} and \cite{X97}.
	\end{proof}
	
	The composition subalgebra of $\cs\cd\widetilde{\ch}(\Lambda)$ is defined to be the image of $\Xi_Q$, which is denoted by $^\imath\widetilde{\cc}(\bfk Q^{\dbl},\swa)$. Then  $\Xi_Q$ induces an isomorphism $\tU|_{v=\sqq}\stackrel{\cong}{\longrightarrow}{}^\imath\widetilde{\cc}(\bfk Q^{\dbl},\swa)$.

	\subsection{Semi-derived Ringel-Hall algebras for $C_n$}
	In this subsection, we shall use the twisted semi-derived Ringel-Hall algebra $\iH(\bfk C_n^{\dbl},\swa)$ to reformulate the results on the twisted extended Ringel-Hall algebra $\widetilde{\ch}(\bfk C_{n})$ of $\rep_{\bfk}^{\rm nil}(C_n)$ that are obtained in \cite{Sch02,Hu05}, by viewing $\widetilde{\ch}(\bfk C_{n})$ as a subalgebra of $\iH(\bfk C_n^{\dbl},\swa)$; see also \cite[Chapter 2]{DDF12}.
	
	For any $m>0$, we define
	\begin{align}
		\label{eq:hr}
		\bh_m=\bh_{0,m}=\frac{[m]_\sqq}{m} \sum_{|\lambda|=m} \bn(\ell(\lambda)-1)\frac{[S_{0}^{(\lambda)}]}{\big|\Aut(S_{0}^{(\lambda)})\big|}\in \widetilde{\ch}(\bfk C_{n})\subset \iH(\bfk C_n^{\dbl},\swa) ;
	\end{align}
	see \cite[(4.1)]{Sch04}.
	It is worth noting that $\bh_{0,m}$ is denoted by $\pi_{1,m}$ in \cite{Hu}.
	
	\begin{lemma}[\text{\cite[Lemma 4.3]{Sch04}}]
		\label{lem:genclassicHallcyclic}
		$\widetilde{\ch}(\bfk C_{n})$ is generated by $\widetilde{\cc}(\bfk C_n)$ and $\bh_m$ for all $m>0$.
	\end{lemma}


	Following Lusztig, let $T_i$ ($0\leq i\leq n-1$) be the braid group actions on $\tU_v(\widehat{\mathfrak{sl}}_n)$; see e.g. \cite[\S 6.2]{DJX12}.
	For $1\leq i\leq n-1$ and $l\in\Z$, let
	\begin{align*}
		x_{i,l}^-=(-1)^{il} T_{\omega_i}^l(F_i),\qquad x_{i,l}^+=(-1)^{il}  T_{\omega_i}^{-l}(E_i),
	\end{align*}
	see \eqref{eq:wj} for the explicit form of $\omega_i$.
	For $m>0$, define
	\begin{align*}
		\psi_{i,m}={(-1)^{im}}(\tK_i')^{-1}[T_{\omega_i}^{-m}(E_i),F_i],
	\end{align*}
	and then define $h_{i,m}$ via
	\begin{align}
		\label{exppsi h}
		1+ \sum_{m\geq 1} (v-v^{-1})\psi_{i,m}u^m &=  \exp\Big((v -v^{-1}) \sum_{m\ge 1}  h_{i,m}u^m\Big).
	\end{align}
	

	Using Lemma \ref{lem:UslnSDH}, for $1\leq i\leq n-1$ and $l\in\Z,m>0$, we define the following root vectors in $\iH(\bfk C_n^{\dbl},\swa)$:
	\begin{align}
		\label{eq:rootxi}
		&\widehat{x}_{i,l}^{+}= (1-q)\Xi_{C_n}(x_{i,l}^{+}), \qquad\qquad \widehat{x}_{i,l}^{-}= (\sqq-\sqq^{-1})\Xi_{C_n}(x_{i,l}^{-}),
		\\
		&\bh_{i,m}=\Xi_{C_n}(h_{i,m}),\qquad\qquad\qquad
		\widehat{\psi}_{i,m}= \Xi_{C_n}(\psi_{i,m}). 
	\end{align}
	
	From \eqref{eq:hax-}--\eqref{eq:hax+} obtained in the proof of Proposition \ref{prop:surjective} below, we know that $\widehat{x}_{i,l}^{+}$, $\widehat{x}_{i,-l}^{-}*[\E_{l\de-\alpha_j}]$, $\widehat{\psi}_{i,m}\in \widetilde{\ch}(\bfk C_{n})$ for any $l,m>0$, (these elements are even in $\widetilde{\cc}(\bfk C_n)$ by \cite{Sch02,Hu05}).

	For $m>0$, we define $\boldsymbol{\pi}_{n,m}\in\widetilde{\ch}(\bfk C_n)$ as
	\begin{align}
		\label{def:pinm}
		\boldsymbol{\pi}_{n,m}=-\sum_{i=0}^{n-1} [n-i]_{\sqq^m}\bh_{i,m}. 
	\end{align}

	\begin{lemma}
		\label{lem:Sch-Hu Theorem}
		The elements $\boldsymbol{\pi}_{n,m}$ ($m>0$) are central in $\widetilde{\ch}(\bfk C_n)$. Moreover, we have
		\begin{align}
			\widetilde{\ch}(\bfk C_n)\cong \widetilde{\cc}(\bfk C_n)\otimes \Q(\sqq)[\boldsymbol{\pi}_{n,m}\mid m>0],
		\end{align}
		where $\Q(\sqq)[\bpi_{n,m}\mid m>0]$ is the polynomial algebra in $\boldsymbol{\pi}_{n,m}$ for $m>0$. In particular,
		$\widetilde{\ch}(\bfk C_n)$ is generated by $[S_i]$ and $\boldsymbol{\pi}_{n,m}$ for $0\leq i\leq n-1$ and $m>0$.
	\end{lemma}

	\begin{proof}
		Recall from \cite[Main Theorem]{Hu} that
		\begin{align}
			\label{him}
			\bh_{i,m}=\pi_{i+1,m}-[2]_{\sqq^m}\pi_{i,m}+\pi_{i-1,m}\quad (1\leq i\leq n-1)
		\end{align}
		in $\widetilde{\ch}(\bfk C_n)$. See \cite[Page 10]{Hu} for the definition of $\pi_{i,m}$.
		Note that $\pi_{0,m}=0$ and $\pi_{1,m}=\bh_{0,m}$. We claim that
		for $1\leq j\leq n$,
		\begin{align}\label{pijm}
			\pi_{j,m}=\sum_{i=0}^{j-1} [j-i]_{\sqq^m}\bh_{i,m}.
		\end{align}
		If the claim holds, then $\bpi_{n,m}=\pi_{n,m}$. So the desired result follows from \cite[Theorem 1]{Hu}.
		
		We prove the claim by induction on $j$.
		If $j=1$, it is obvious. Now assume \eqref{pijm} holds for any $1\leq j'<j$. Then from \eqref{him} we have
		\begin{align*}
			\pi_{j,m}=&\bh_{j-1,m}+[2]_{\sqq^m}\pi_{j-1,m}-\pi_{j-2,m}\\
			=&\bh_{j-1,m}+[2]_{\sqq^m}\sum_{i=0}^{j-2} [j-1-i]_{\sqq^m}\bh_{i,m}-\sum_{i=0}^{j-3} [j-2-i]_{\sqq^m}\bh_{i,m}\\
			=&\bh_{j-1,m}+[2]_{\sqq^m}\bh_{j-2,m}+\sum_{i=0}^{j-3}  \big([2]_{\sqq^m}\cdot [j-1-i]_{\sqq^m}-[j-2-i]_{\sqq^m}\big)\bh_{i,m}\\
			=&\bh_{j-1,m}+[2]_{\sqq^m}\bh_{j-2,m}+\sum_{i=0}^{j-3}  [j-i]_{\sqq^m}\bh_{i,m}\\
			=&\sum_{i=0}^{j-1}  [j-i]_{\sqq^m}\bh_{i,m},
		\end{align*}
		where the last second equality follows from the general fact $[2]\cdot [i]=[i-1]+[i+1]$.
		This proves the claim, and then the result follows.
	\end{proof}

	\subsection{$\imath$Quantum loop algebra $\tUiD_v(\widehat{\mathfrak{gl}}_n)$}
	Recall the presentation of $\tUiD_v(\widehat{\mathfrak{sl}}_n)$ in Definition \ref{def:iDR}.
	Let us give the definition of $\imath$quantum loop algebras $\tUiD_v(\widehat{\mathfrak{gl}}_n)$ for $n>1$.
	\begin{definition}
		\label{def:loopgln}
		The $\imath$quantum loop algebra $\tUiD_v(\widehat{\mathfrak{gl}}_n)$ is the $\Q(v)$-algebra generated by $\tUiD_v(\widehat{\mathfrak{sl}}_n)$ and $H_{0,m}$ for $m>0$, subject to the following relations for $0\leq i\leq n-1$, $r,m>0$, $l\in\Z$:
		\begin{align}
			\label{eq:HHcom}
			&\K_i \text{ and } C \text{ are central,}\qquad
			[H_{0,m},H_{i,r}]=0,
			\\
			\label{HBcom}
			&[H_{0,m},\y_{j,l}]=\begin{cases}-\frac{[m]}{m} \y_{1,l+m}+\frac{[m]}{m} \y_{1,l-m}C^m, & \text{ if }j=1,
				\\
				0,& \text{ otherwise}.
			\end{cases}
		\end{align}
	\end{definition}
	

	Inspired by \cite{Hu} (see also \cite[Chapter 2]{DDF12}), for any $m>0$, define
	\begin{align}
		\theta_m:=\sum_{i=0}^{n-1} [n-i]_{v^m}H_{i,m}\in\tUiD_v(\widehat{\mathfrak{gl}}_n).
	\end{align}
	
	\begin{lemma}
		\label{lem:central}
		For any $m>0$, $\theta_m$ is a central element of $\tUiD_v(\widehat{\mathfrak{gl}}_n)$.
	\end{lemma}
	
	\begin{proof}
		It is enough to prove that $[\theta_m,B_{j,l}]=0$ for any $1\leq j\leq n-1$ and $l\in\Z$.
		In fact, if $j>1$, by $[H_{0,m},B_{j,l}]=0$, we have
		\begin{align*}
			[\theta_m,B_{j,l}]=&  \sum_{i=0}^{n-1} [n-i]_{v^m}[H_{i,m},B_{j,l}]
			\\
			=&\sum_{i=1}^{n-1} [n-i]_{v^m}\frac{[mc_{ij}]}{m} \cdot\Big( \y_{j,l+m}-\y_{j,l-m}C^m\Big)
			\\
			=&0.
		\end{align*}
		Here we use
		\begin{align*}
			\sum_{i=1}^{n-1}[n-i]_{v^m}[mc_{ij}] =& -[n-j+1]_{v^m} [m]+[n-j]_{v^m} [2m]-[n-j-1]_{v^m} [m]
			=0.
		\end{align*}
		
		If $j=1$, the proof is similar, hence omitted here.
	\end{proof}

	\subsection{Isomorphism between  $\tUiD_v(\widehat{\mathfrak{gl}}_n)$ and $\iH(\bfk C_n)$}
	
	Recall that ${}^{\imath}\widetilde{\ch}(\bfk C_{n})$ is the $\imath$Hall algebra of $\rep^{\rm nil}_\bfk(C_n)$.
	Inspired by \eqref{def:Hxm}, for $m>0$ we define
	\begin{align}
		\label{def:H0m}
		\widehat{H}_{0,m}:=\frac{[m]_{\sqq}}{m} \sum_{|\lambda|=m} \bn(\ell(\lambda)-1)\frac{[S_{0}^{(\lambda)}]}{\big|\Aut(S_{0}^{(\lambda)})\big|}-\delta_{m, ev} \sqq^{-\frac{m}{2}} \frac{[m/2]_{\sqq}}{m} [K_{\frac{m}{2}\de}]\in {}^{\imath}\widetilde{\ch}(\bfk C_{n}).
	\end{align}
	
	Recall $\Omega_{C_n}: \tUiD_\sqq(\widehat{\mathfrak{sl}}_n)\rightarrow \iH(\bfk C_n)$ in \eqref{the map Psi A}. In the following lemma, we extend $\Omega_{C_n}$ to be defined on $\tUiD_v(\widehat{\mathfrak{gl}}_n)$.
	\begin{lemma}
		\label{lem:morCn}
		There exists an algebra morphism $\Omega_{C_n}: \tUiD_\sqq(\widehat{\mathfrak{gl}}_n)\rightarrow \iH(\bfk C_n)$, which sends
		\begin{align}
			&\K_j\mapsto [K_{S_j}],\qquad  B_{j,0}\mapsto \frac{-1}{q-1}[S_j], \qquad \forall \,\,0\leq j\leq n-1,\\
			&C\mapsto [K_\de],\qquad
			H_{j,m}\mapsto \widehat{H}_{j,m}, \text{ for }0\leq j\leq n-1 \text{ and } m>0.
		\end{align}
	\end{lemma}
	
	\begin{proof}
		By the definition of root vectors, it is enough to check that $\Omega_{C_n}$ preserves \eqref{eq:HHcom} and \eqref{HBcom}.
		Assume that $\X_\bfk$ is a weighted projective line with $p_1=n$. Then it implies that
		$\iH(\bfk C_n)$ is a subalgebra of $\iH(\X_\bfk)$. By \cite[Propositions 7.8, 8.6]{LR21}, we have
		$$[\haH_{\ast,r},\haH_{i,m}]=0, \text{ for }r,m\geq0,\;\; 1\leq i\leq n-1.$$
		Observe that there are no non-zero homomorphisms and extensions between distinct tubes.	Using \eqref{eq:Hdecomptox}, we have
		$$[\haH_{0,r},\haH_{i,m}]=0, \text{ for }r,m\geq0,\;\; 0\leq i\leq n-1.$$
		
		Using the fact that $\iH(\bfk C_1)$ is commutative (see \cite[Proposition 4.5]{LRW20a}), together with the natural embedding $F_{\X,\P^1}$, we have 
		$$[\haH_{0,r},\haH_{0,m}]=0, \text{ for }r,m\geq0.$$
		
		Similarly, for $m>0,l\in\Z$, we have
		$$[\haH_{0,m},\haB_{j,l}]=[\haH_{\star,m},\haB_{[1,j],l}],$$
		and then the desired formula follows from \cite[Propositions 7.7, 8.5]{LR21}.
	\end{proof}

	Let
	\begin{align}
		\label{def:hathetam}
		\hath_m:=\Omega_{C_n}(\theta_m)=\sum_{i=0}^{n-1} [n-i]_{\sqq^m}\haH_{i,m}.
	\end{align}
	
	\begin{proposition}
		\label{prop:surjective}
		The  algebra morphism $\Omega_{C_n}: \tUiD_\sqq(\widehat{\mathfrak{gl}}_n)\rightarrow \iH(\bfk C_n)$ defined in Lemma \ref{lem:morCn} is surjective. Moreover, $\hath_m$ is central for $m>0$.
	\end{proposition}

	\begin{proof}
		Recall the $\imath$Hall basis obtained in Proposition \ref{prop:hallbasis}~(3). For any $\gamma\in K_0(\rep^{\rm nil}_\bfk(C_n))$, denote by $\iH(\bfk C_n)_{\leq\gamma}$ (respectively, $\iH(\bfk C_n)_{\gamma}$) the subspace of $\iH(\bfk C_n)$ spanned by elements from this basis for which $\widehat{M}\leq \gamma$ (respectively, $\widehat{M}= \gamma$) in $K_0(\rep^{\rm nil}_\bfk(C_n))$.
		Then $\iH(\bfk C_n)$ is a $K_0(\rep^{\rm nil}_\bfk(C_n))$-graded linear space, i.e.,
		\begin{align}
			\label{eqn: graded linear MH}
			\iH(\bfk C_n)=\bigoplus_{\gamma\in K_0(\rep^{\rm nil}_\bfk(C_n))} \iH(\bfk C_n)_\gamma.
		\end{align}
		By \cite[Lemma 5.4]{LW19a}, $\iH(\bfk C_n)$ is a filtered algebra with its associated graded algebra $\iH(\bfk C_n)^{\gr}$ isomorphic to $\widetilde{\ch}(\bfk C_n)\otimes \Q(\sqq)\big[[K_i]^{\pm1}\mid 0\leq i<n\big]$. Using Lemma \ref{lem:Sch-Hu Theorem}, we know that $\widetilde{\ch}(\bfk C_n)$  is generated by
		$[S_i]$, $\bpi_{n,m}$ for $0\leq i< n$ and $m>0$. 
		
		Recall $\bpi_{n,m}$ defined in \eqref{def:pinm}. View $\bpi_{n,m}$ in $\iH(\bfk C_n)$.
		We claim that for any $m>0$,
		\begin{align}
			\label{eq:theta pi}
			\hath_m-\bpi_{n,m}\in \iH(\bfk C_n)_{<m\de}.
		\end{align}
		If the claim holds, then the first statement follows by a standard argument of filtered algebras, and then the second statement follows from the first one and Lemma \ref{lem:central}.
		
		Let us prove the claim in the following. Compare \eqref{def:pinm} with \eqref{def:hathetam}, and \eqref{eq:hr} with \eqref{def:H0m}. It is enough to prove that
		\begin{align}
			\haH_{j,m}-\bh_{j,m}\in \iH(\bfk C_n)_{<m\de},\;\;\;\forall 1\leq j\leq n-1.
		\end{align}

		Following \cite{LW19b} (see also \cite{Rin96,SV99,XY01}), we can interpret $T_i$ in the twisted semi-derived Ringel-Hall algebras
		by using the BGP type reflection functors. Let $Q$ be a quiver with its underlying graph the same as $C_n$. Then there exists an isomorphism $\texttt{FT}_{C_n,Q}:{}^\imath\widetilde{\cc}(\bfk Q^{\dbl},\swa)\stackrel{\cong}{\longrightarrow} {}^\imath\widetilde{\cc}(\bfk C_n^{\dbl},\swa)$ 
		such that $\texttt{FT}_{C_n,Q}([S_i])=[S_i]$, $\texttt{FT}_{C_n,Q}([S_{i'}])=[S_{i'}]$, $\texttt{FT}_{C_n,Q}([\E_i])=[\E_i]$, $\texttt{FT}_{C_n,Q}([\E_{i'}])=[\E_{i'}]$ for any $0\leq i\leq n-1$; cf. Lemma \ref{lem:UslnSDH}. Recall the quiver $Q(j)$  defined in \eqref{quiver:affineA}. In order to avoid confusion, we denote ${}^j\Lambda:=\bfk Q(j)\otimes_\bfk R_2$. Similarly to Proposition \ref{rem:rootvectors}, for any $l>0$, we have
		\begin{align}
			\label{eq:hax-}
			\widehat{x}_{j,l}^-=& (-1)^{jl}[\E_{l\de-\alpha_j}]^{-1}*\texttt{FT}_{C_n,Q(j)}([M(l\de-\alpha_j)]),
			\\
			\label{eq:hax+}
			\widehat{x}_{j,l}^+= &(-1)^{jl}\texttt{FT}_{C_n,Q(j)}([M(l\de+\alpha_j)]).
		\end{align}
		
		By using \cite[Lemma 4.3]{LW19a}, it follows that
		\begin{align*}
			\widehat{\psi}_{j,l}=&\Xi_{C_n}(\psi_{j,l})
			\\
			=&{(-1)^{jl}}\Xi_{C_n}\big((\tK_j')^{-1}[T_{\omega_j}^{-m}(E_j),F_j]\big)
			\\
			=&\frac{{(-1)^{jl}}}{1-q}\cdot\frac{\sqq}{q-1}[\E_{j}]^{-1}*\texttt{FT}_{C_n,Q(j)}\Big(\big[[M(l\de+\alpha_j)], [S_{j'}] \big]\Big)
			\\
			=&\frac{{(-1)^{jl}}\sqq}{(q-1)^2}[\E_{j}]^{-1}*\texttt{FT}_{C_n,Q(j)}\Big(\sum_{0\neq g: S_j \rightarrow M(l\de+\alpha_j)} [\E_j\oplus\coker g]\Big).
		\end{align*}
		Here the last equality follows by noting that the class of a non-split short exact sequence $0\rightarrow S_{j'}\rightarrow N\rightarrow M(l\de+\alpha_j)\rightarrow0$ in $\Ext^1_{{}^j\Lambda}(M(l\de+\alpha_j),S_{j'})$ is uniquely determined by a nonzero morphism $g:S_j \rightarrow M(l\de+\alpha_j)$. In this case, $N$ admits a submodule $\E_j$, and then $[N]=[\E_j\oplus\coker g]$ in $\cs\cd\widetilde{\ch}({}^j\Lambda)$. Then
		\begin{align*}
			\widehat{\psi}_{j,l}=&\frac{{(-1)^{jl}}\sqq}{(q-1)^2}\sqq^{-\langle S_j,M(l\de)\rangle_{Q(j)}}q^{\langle \E_j,M(l\de)\rangle_{{}^j\Lambda} }[\E_{j}]^{-1}*\texttt{FT}_{C_n,Q(j)}\Big(\sum_{0\neq g: S_j \rightarrow M(l\de+\alpha_j)} [\E_j]*[\coker g]\Big)
			\\
			=&\frac{{(-1)^{jl}}\sqq}{(q-1)^2}\sqq^{-\langle S_j,M(l\de)\rangle_{Q(j)}}\texttt{FT}_{C_n,Q(j)}\Big(\sum_{0\neq g: S_j \rightarrow M(l\de+\alpha_j)} [\coker g]\Big)
			\\
			=&\frac{{(-1)^{jl}}}{(q-1)^2\sqq^{l-1}}\texttt{FT}_{C_n,Q(j)}\Big(\sum_{0\neq g: S_j \rightarrow M(l\de+\alpha_j)} [\coker g]\Big).
		\end{align*}
		Here $\langle \E_j,M(l\de)\rangle_{{}^j\Lambda} =0$ by using \cite[Proposition 2.4]{LP16}; and $\langle S_j,M(l\de)\rangle_{Q(j)}=l$.
		
		View $\widetilde{\ch}(\bfk C_n)$ as the subalgebra of $\iH(\bfk C_n)^{\gr}$. For any element $Y\in \iH(\bfk C_n)$, denote by $\gr(Y)$ the corresponding element in $\iH(\bfk C_n)^{\gr}$.
		By Proposition \ref{rem:rootvectors}, we have 
		$\gr(\haTh_{j,l})=\widehat{\psi}_{j,l}$ for any $l>0$.
		As a consequence, $\gr(\haH_{j,l})= \bh_{j,l}$  by comparing the counterparts of \eqref{exppsi h} and \eqref{exp h} in Hall algebras. Therefore, the claim follows.
	\end{proof}


	From Proposition \ref{prop:surjective}, we have the following corollary.
	
	\begin{corollary}
		\label{cor:iHallCn}
		We have
		\begin{align}
			\label{eq:Ugln}
			\iH(\bfk C_n)\cong {}^\imath \widetilde{\cc}(\bfk C_n)\otimes \Q(\sqq)[\hath_m\mid m>0 ].
		\end{align}
	\end{corollary}
	
	\begin{proof}
		Note that $\iH(\bfk C_n)^{\gr}\cong \widetilde{\ch}(\bfk C_n)\otimes \Q(\sqq)\big[[K_i]^{\pm1}\mid 0\leq i<n\big]$.
		Comparing Proposition \ref{prop:surjective} with Lemma \ref{lem:Sch-Hu Theorem}, and using \eqref{eq:theta pi}, the desired result follows by a standard argument of filtered algebras.
	\end{proof}
	
	\begin{theorem}
		\label{prop:iso}
		For any $n\geq2$, the  algebra morphism $\Omega_{C_n}: \tUiD_\sqq(\widehat{\mathfrak{gl}}_n)\rightarrow \iH(\bfk C_n)$ obtained in Lemma \ref{lem:morCn} is an isomorphism.
	\end{theorem}

	\begin{proof}
		Applying \eqref{eq:isocomp} to the cyclic quiver $C_n$, we have 
		${}^\imath \widetilde{\cc}\cong \tUi_\sqq(\widehat{\mathfrak{sl}}_n)$. 
		Using \eqref{eq:Ugln}, we define $\Psi_{C_n}: \iH(\bfk C_n) \rightarrow \tUiD_\sqq(\widehat{\mathfrak{gl}}_n)$ by sending
		\begin{align*}
			& [K_{S_j}]\mapsto\K_j,\qquad [K_\de]\mapsto C,\qquad \frac{-1}{q-1}[S_j]\mapsto B_{j,0}, \\
			&
			\widehat{H}_{j,m}\mapsto H_{j,m},\qquad \hath_m\mapsto \theta_m,
		\end{align*}
		for $1\leq j\leq n-1 \text{ and } m>0$.
		Using Lemma \ref{lem:central} and definitions of root vectors in $\iH(\bfk C_n)$,  one can check that $\Psi_{C_n}$ is well defined.
		
		It is easy to see that $\tUiD_v(\widehat{\mathfrak{gl}}_n)$ is generated by $\tUi_v(\widehat{\mathfrak{sl}}_n)$ and $\theta_{m}$ ($m>0$). So one can check that $\Psi_{C_n}$ is the inverse map of $\Omega_{C_n}: \tUiD_\sqq(\widehat{\mathfrak{gl}}_n)\rightarrow \iH(\bfk C_n)$, and then $\Omega_{C_n}$ is an isomorphism. 
	\end{proof}
	
	Using Theorem \ref{prop:iso} and the definition of  $\tUiD_v(\widehat{\mathfrak{gl}}_n)$, we have the following corollary.
	\begin{corollary}
		\label{cor:present HCn}
		$\iH(\bfk C_n)$ is isomorphic to the $\Q(\sqq)$-algebra generated by ${}^\imath \widetilde{\cc}(\bfk C_n)$ and $\widehat{H}_{0,m}$ ($m>0$) subject to the following relations for $0\leq i\leq n-1$, $r,m>0$, $l\in\Z$:
		\begin{align}
			&[K_{S_i}] \text{ are central}, \qquad [\haH_{0,m},\widehat{H}_{i,r}]=0,
			\\
			&[\haH_{0,m},\haB_{j,l}]=\begin{cases}\frac{-[m]_\sqq}{m} \haB_{1,l+m}+\frac{[m]_\sqq}{m} \haB_{1,l-m}*[K_\de]^m,& \text{ if }j=1,
				\\
				0&\text{ otherwise}.
			\end{cases}
		\end{align}
	\end{corollary}

	\subsection{$\imath$Quantum group $\tUi_v(\widehat{\mathfrak{gl}}_n)$}
	
	\begin{definition}
		The $\imath$quantum group of affine $\mathfrak{gl}_n$ is the $\Q(v)$-algebra $\tUi_v(\widehat{\mathfrak{gl}}_n)$ generated by $B_i,\K_i,\K_i',c_m$, for $0\leq i\leq n-1$, $m>0$ with the relations given by \eqref{eq:KK}--\eqref{eq:S3}, and
		\begin{align}
			& c_{m} \text{ is central, for } m>0.
		\end{align}
	\end{definition}
	By definition, we have that $\tUi_v(\widehat{\mathfrak{sl}}_n)$ is naturally a subalgebra of $\tUi_v(\widehat{\mathfrak{gl}}_n)$. Moreover,
	\begin{align}
		\tUi_v(\widehat{\mathfrak{gl}}_n)\cong \tUi_v(\widehat{\mathfrak{sl}}_n)\otimes_{\Q(v)} \Q(v)[c_{m}\mid m>0].
	\end{align}

	\begin{remark}
		\label{rem:iso 2presentations}
		In a separate publication, we shall prove that $\tUi_v(\widehat{\mathfrak{gl}}_n)$ is isomorphic to $\tUiD_v(\widehat{\mathfrak{gl}}_n)$, which justifies the terminologies; cf. Theorem \ref{prop:iso} and Corollary \ref{cor:iHallCn}. In particular, the relation between the central elements $\theta_m$ and $c_m$ shall be described clearly; cf.
		\cite{Hu}.
	\end{remark}

	\begin{remark}
		For any $m>0$, define
		\begin{align}
			\label{def:cm}
			\widehat{c}_{m}=(-1)^m\sqq^{-2nm}\sum_{[M]} (-1)^{\dim \End(M)}[M]\in \iH(\bfk C_n),
		\end{align}
		where the sum is taken over all $[M]\in\Iso(\rep_\bfk^{\rm nil}(C_n))$ such that $\widehat{M}=m\de$ and ${\rm soc} M\subseteq S_1\oplus \cdots \oplus S_n$.
		We also define $\widehat{c}_0=1$ by convention. 
		In a separate publication, we shall prove that $\widehat{c}_m (m>0)$ are central elements, and
		\begin{align}
			\iH(\bfk C_n)\cong {}^\imath \widetilde{\cc}(\bfk C_n)\otimes \Q(\sqq)[\widehat{c}_m\mid m>0].
		\end{align}
		In particular, we shall prove that $\tUi_\sqq(\widehat{\mathfrak{gl}}_n)\cong \iH(\bfk C_n)$, which maps $c_m$ to $\widehat{c}_m$. (This isomorphism is compatible with the ones given in Theorem \ref{prop:iso} and Remark \ref{rem:iso 2presentations}.)
	\end{remark}

	\section{Injectivity for finite and affine type cases}
	\label{sec:injectivity}
	
	In this section, we shall prove that the algebra homomorphism $\Omega: \tUiD_{|v=\sqq}\rightarrow \iH(\X_\bfk)$ defined in \S\ref{subsec:homo} is injective if the star-shaped graph $\Gamma$ is of finite or affine type.
	
	%
	%
	\subsection{$\imath$Composition subalgebra}
	
	Recall the definition of $\tCMHX$ in \S\ref{subsec:homo}.
	In this subsection, we give some general results on the composition subalgebras $\tCMHX$ for arbitrary $\X_\bfk$.
	
	\begin{lemma}
		\label{lem:linebundle}
		The algebra $\tCMHX$ contains the classes $[\co(\vec{x})]$ of all line bundles in $\coh(\X)$.
	\end{lemma}
	
	\begin{proof}
		Assume that $\vec{x}=\sum_{i=1}^\bt l_i\vec{x}_i+l\vec{c}$ in normal form. We prove 	$[\co(\vec{x})]\in \tCMHX$ by induction on $\sum_{i=1}^\bt l_i$. In fact, if $\sum_{i=1}^\bt l_i=0$, then $[\co(\vec{x})]=[\co(l\vec{c})]\in \tCMHX$ by definition.
		
		If $\sum_{i=1}^\bt l_i>0$, then there exists some $l_i>0$. Consider the following exact sequence
		\begin{align*}
			0\rightarrow \co(\vec{x}-\vec{x}_i)\rightarrow \co(\vec{x})\rightarrow S_{i,l_i}\rightarrow0.
		\end{align*}
		We obtain that
		\begin{align*}
			[\co(\vec{x})]=\frac{1}{\sqq-\sqq^{-1}}\big[[S_{i,l_i}],[\co(\vec{x}-\vec{x}_i)]\big]_{\sqq^{-1}}\in\tCMHX
		\end{align*}
		by
		the inductive assumption. This finishes the proof.
	\end{proof}

	\begin{lemma}
		\label{lem:genlinebundle}
		$\tCMHX$ is generated by $[\co(\vec{x})]$ $(\vec{x}\in L(\bp))$ and $[K_\beta]$  $(\beta\in K_0(\coh(\X_\bfk)))$.
	\end{lemma}
	
	\begin{proof}
		Let $\bH$ be the subalgebra of  $\tCMHX$ generated by $[\co(\vec{x})]$ $(\vec{x}\in L(\bp))$ and $[K_\beta]$  $(\beta\in K_0(\coh(\X_\bfk)))$.
		By definition, it is enough to check that all the isoclasses of simple sheaves $[S_{ij}]$ for $1\leq i\leq \bt$ and $0\leq j\leq p_i-1$ are in $\bH$.
		
		For $[S_{ij}]$, we have a short exact sequence
		$$0\longrightarrow \co((j-1)\vec{x}_i)\stackrel{x_i}{\longrightarrow} \co(j\vec{x}_i)\longrightarrow S_{ij}\rightarrow0.$$
		Then
		\begin{align*}
			\big[[\co((j-1)\vec{x}_i)],[\co(j\vec{x}_i)]\big]_{\sqq^{-1}}=&(\sqq-\sqq^{-1})[S_{ij}]*[K_{\co((j-1)\vec{x}_i)}].
		\end{align*}
		So $[S_{ij}]\in \bH$. The proof is completed.
	\end{proof}

	\subsection{A Drinfeld type presentation of $\widetilde{\cc}(\X_\bfk)$}

	Let $\Gamma=\T_{p_1,\dots,p_\bt}$ be the star-shaped graph defined in \eqref{star-shaped}. Let $\widetilde{\ch}(\X_\bfk)$ be the (twisted) Ringel-Hall algebra of $\coh(\X_\bfk)$; cf. \eqref{eq:RHallmult}.  Then $\widetilde{\ch}(\X_\bfk)$ is a $K_0(\coh(\X_\bfk))$-graded algebra.
	Note that $\rep_\bfk^{\rm nil}(C_{p_i})\simeq \scrt_{\bla_i}$, which induces an algebra embedding
	$\iota_i:\widetilde{\ch}(\bfk C_{p_i})\rightarrow \widetilde{\ch}(\X_\bfk)$.
	Inspired by \eqref{eq:rootxi}, define
	\begin{align}
		\label{eq:Eijk}
		\widehat{x}^+_{[i,j],k}:=(1-q) \iota_i(\widehat{x}^+_{j,k})\in \widetilde{\ch}(\X_\bfk),\quad \forall k\geq0.
	\end{align}

	We shall recall the algebra $\U_\sqq(\hn)$ defined in \cite{Sch04}, which can be viewed as a ``positive part'' of the quantum loop algebra associated to $\Gamma$.
	
	Denote by $\phi_i: \U_\sqq^+(\widehat{\mathfrak{sl}}_{p_i})\rightarrow \widetilde{\ch}(\bfk C_{p_i})$ the algebra embedding obtained in \eqref{eq:Xi}. We denote by $E_0^{(i)}$, $E_{\phi_i(1)}$, $\cdots$, $E_{\phi_i(p_i-1)}$ the corresponding elements $E_0,E_1,\cdots, E_{p_i-1}$ under the morphism $\phi_i$, for $1\leq i\leq \bt$.
	Inspired by \eqref{eq:hr}, we define
	\begin{align}
		\bh_{[i,0],r}=\frac{[r]_\sqq}{r} \sum_{|\lambda|=r} \bn(\ell(\lambda)-1)\frac{[S_{i,0}^{(\lambda)}]}{\big|\Aut(S_{i,0}^{(\lambda)})\big|}\in \widetilde{\ch}(\bfk C_{p_i})\subset \widetilde{\ch}(\X_\bfk).
	\end{align}
	
	Let $\U_\sqq(\hn)$ be the $\Q(\sqq)$-algebra generated by $\widetilde{\ch}(\bfk C_{p_i})$ for $1\leq i\leq \bt$, $\bh_{\star,r}$ for $r>0$ and $E_{\star,k}$ for $k\in\Z$, subject to a set of relations defined in \cite[\S 4.4]{Sch04}.

	Let $\tCHX$ be the composition subalgebra of $\widetilde{\ch}(\X_\bfk)$ defined in \cite{Sch04}, i.e.,
	$\tCHX$ is generated by
	$[\co(l\vec{c})]$, $[S_{ij}]$, $T_r$ for $l\in\Z$, $r>0$, $1\leq i\leq \bt$, $0\leq j\leq p_i-1$.
	Here
	\begin{align*}
		T_{r}
		:=&\sum_{x,d_x|r} \frac{[r]_\sqq}{r} d_x \sum_{|\lambda|=\frac{r}{{d_x}}} \bn_x(\ell(\lambda)-1)\frac{[\mathbb{F}_{\X,\P^1}\big(S_x^{(\lambda)}\big)]}{\big|\Aut(S_x^{(\lambda)})\big|};
	\end{align*}
	see \eqref{the embedding functor F} for the definition of $\mathbb{F}_{\X,\P^1}$.
	It is a fact that $\widetilde{\cc}(\bfk C_{p_i})$ is a subalgebra of $\tCHX$ by using $\iota_i$.
	So $\widehat{x}^+_{[i,j],k}\in \tCHX$; see \eqref{eq:Eijk}.

	\begin{theorem}
		[\cite{Sch04}]
		\label{thm:Sch}
		There is an algebra epimorphism $\Phi: \U_\sqq(\widehat{\bn})\rightarrow \tCHX$ uniquely determined by
		\begin{align*}
			E_{\star,k}\mapsto \frac{-1}{q-1}[\co(k\vec{c})],
			\qquad E_{\phi_i(j)}\mapsto \frac{-1}{q-1}[S_{ij}], \qquad E_0^{(i)}\mapsto \frac{-1}{q-1} [S_{i,0}],  \qquad
			\bh_{\star,r}\mapsto T_r.
		\end{align*}
		Moreover, $\Phi$ is isomorphic if $\Gamma$ is of finite or affine type.
	\end{theorem}

	We assume that $\Gamma$ is of finite or affine type in the following. Using the presentation of $\U_\sqq(\widehat{\bn})$ defined in \cite[\S 4.4]{Sch04} and the isomorphism $\Phi: \U_\sqq(\widehat{\bn})\rightarrow \tCHX$ in Theorem \ref{thm:Sch}, we obtain a presentation of $\tCHX$.
	
	\begin{corollary}
		\label{cor:presentationU+}
		If $\Gamma$ is of finite or affine type, then $\tCHX$ is isomorphic to the algebra $\widetilde{\bC}_{\X_\bfk}$ over $\Q(\sqq)$ generated by $\widetilde{\ch}(\bfk C_{p_i})$, $[\co(l\vec{c})]$, $T_r$ for $1\leq i\leq \bt$, $l\in\Z$, $r>0$, subject to the following relations precisely:
		\begin{align}
			\label{eq:U+Dr1}
			&T_r=\bh_{[i,0],r},  \qquad [\widetilde{\cc}(\bfk C_{p_i}),\widetilde{\cc}(\bfk C_{p_j}) ]=0, \text{ for }1\leq i\neq j\leq \bt;
			\\
			\label{eq:U+Dr2}
			&\big[T_{r}, [\co(l\vec{c})]\big]=\frac{[2r]_\sqq}{r}[\co((l+r)\vec{c})],
			\\
			\label{eq:U+Dr3}
			&\big[[\co(k\vec{c})], [\co((l+1)\vec{c})]\big]_{\sqq^{-2}}  -\sqq^{-2} \big[[\co((k+1)\vec{c})], [\co(l\vec{c})]\big]_{\sqq^{2}}=0, \text{ for }k,l\in\Z,
			\\
			\label{eq:U+Dr4}
			&\big[[\co(l\vec{c})],[S_{i,0}]\big]_{\sqq^{-1}}=0, \qquad \big[[\co(l\vec{c})],[S_{ij}]\big]=0, \text{ if }j>1, \\
			\label{eq:U+Dr5}
			&\Big[[S_{i,0}], \big[ [S_{i,1}],[\co(l\vec{c})]\big]_{\sqq^{-1}} \Big]=0,
			\text{ if } p_i>2,
			\\
			\label{eq:U+Dr8}
			&\big[ [\co((l+1)\vec{c})],\widehat{x}^+_{[i,1],r}\big]_{\sqq^{-1}}-\sqq^{-1} \big[[\co(l\vec{c})],\widehat{x}^+_{[i,1],r+1}\big]_\sqq=0,
			\\
			\label{eq:U+Dr6}
			&\widehat{\SS}(t_1,t_2\mid l;[i,1],\star )=0,
			\text{ for }t_1,t_2>0,\,\, l\in\Z,
			\\
			\label{eq:U+Dr7}
			&\widehat{\SS}(l_1,l_2\mid r;\star,[i,1])=0,
			\text{ for } l_1,l_2\in\Z, \,\, r>0.
		\end{align}
	\end{corollary}
	




	\subsection{Filtered algebra structure of $^\imath\widetilde{\cc}(\X_\bfk)$}
	Recall from Theorem \ref{prop:hallbasis} that $\iH(\X_\bfk)$ has an $\imath$Hall basis given by $[M]* [K_\alpha]$, where $[M]\in\Iso(\coh(\X_\bfk))$, and $\alpha\in K_0(\coh(\X_\bfk))$.
	For any $\gamma\in K_0(\coh(\X_\bfk))$, denote by $\iH(\X_\bfk)_{\leq\gamma}$ (respectively, $\iH(\X_\bfk)_{\gamma}$) the subspace of $\iH(\X_\bfk)$ spanned by elements from this basis for which $\widehat{M}\leq \gamma$ (respectively, $\widehat{M}= \gamma$) in $K_0(\coh(\X_\bfk))$.
	Then $\iH(\X_\bfk)$ is a $K_0(\coh(\X_\bfk))$-graded linear space, i.e.,
	\begin{align}
		\label{eqn: graded linear MH}
		\iH(\X_\bfk)=\bigoplus_{\gamma\in K_0(\coh(\X_\bfk))} \iH(\X_\bfk)_\gamma.
	\end{align}
	Similar to \cite[Lemma 5.3]{LW19a}, there is a filtered algebra structure on $\tMHX$, and we denote the associated graded algebra by
	\[
	\tMHX^{\gr} = \bigoplus_{\gamma \in K_0(\coh(\X_\bfk))}\tMHX_{\gamma}^{\gr}.
	\]
	It is natural to view the quantum torus $\widetilde{\ct}(\X_\bfk):=\widetilde{\ct}(\coh(\X_\bfk))$ as a subalgebra of $\tMHX^{\gr}$. Then $\tMHX^{\gr}$ is also a $\widetilde{\ct}(\X_\bfk)$-bimodule.

	Just as in \cite[Lemma 5.4 (ii)]{LW19a}, the linear map
	\begin{align}
		\label{def:morph}
		\varphi: \widetilde{\ch}(\X_\bfk)\longrightarrow \tMHX^{\gr},\quad
		\varphi([M])=[M],\;  \forall M\in \coh(\X_\bfk),
	\end{align}
	is an algebra embedding. Furthermore,
	$\varphi$ induces an isomorphism of algebras
	$$\widetilde{\ch}(\X_\bfk)\otimes_{\Q(\sqq)} \Q(\sqq)\big[[K_i]^{\pm1}\mid i\in\I\big]\stackrel{\cong}{\longrightarrow} \tMHX^{\gr},$$
	which is also denoted by $\varphi$.

	Recall that $\tCMHX$ is the composition subalgebra of $\tMHX$. Then $\tCMHX$ is also a filtered algebra, we denote by $\tCMHX^{\gr}$ its associated graded algebra.
	For any element $Y\in \tCMHX$, denote by $\gr(Y)$ the corresponding element in $\tCMHX^{\gr}$.
	\begin{lemma}
		\label{lem:phirestrction}
		The morphism $\varphi$ restricts to an embedding of algebras
		\begin{align}
			\label{eq:embedd}
			\varphi: \tCHX\longrightarrow \tCMHX^{\gr},
		\end{align}
		which induces an isomorphism
		\begin{align}
			\label{eq:isoHiH}
			\tCHX\otimes\Q(\sqq)\big[[K_i]^{\pm1}\mid i\in\I\big] \cong  \tCMHX^{\gr}.
		\end{align}
	\end{lemma}
	
	\begin{proof}
		For \eqref{eq:embedd}, by definitions, using \eqref{def:morph}, it is enough to prove that $[S_{i,0}]\in \tCMHX$, which follows by the proof of Lemma \ref{lem:genlinebundle}.
		
		Let us prove \eqref{eq:isoHiH}. By Corollary \ref{cor:epimorphism} and Definition \ref{def:iDR}, we know  $\tCMHX$ is generated by $[K_\beta]$ ($\beta\in K_0(\coh(X))$) and
		\begin{align}
			[\co(l\vec{c})],\,{\widehat{\Theta}_{\star,r}},\,\haB_{[i,j],l},\,\widehat{\Theta}_{[i,j],r},\; \,\,\forall
			[i,j]\in\II-\{\star\},l\in\Z,r>0.
		\end{align}
		By using Corollary \ref{cor:epimorphism} again, it follows from \cite[Chapter 7, Corollary 6.14]{MR01} that $\tCMHX^{\gr}$ is generated by
		\begin{align}
			\label{eq:generators}
			\gr([\co(l\vec{c})]),\,\gr({\widehat{\Theta}_{\star,r}}),\,\gr(\haB_{[i,j],l}),\,\gr(\widehat{\Theta}_{[i,j],r}),\; \,\,\forall
			[i,j]\in\II-\{\star\},l\in\Z,r>0.
		\end{align}
		
		Obviously,  $\gr([\co(l\vec{c})]),\,\gr({\widehat{\Theta}_{\star,r}})\in\varphi(\tCHX)$. Note that ${}^\imath\widetilde{\cc}(\bfk C_{p_i})$ can be viewed as a subalgebra of $\tCMHX$ via the embedding $\iota_i$ defined in \eqref{eq:embeddingx}. From the proof of Proposition \ref{prop:surjective},  we have
		$\gr(\haB_{[i,j],l}),\,\gr(\haTh_{[i,j],r})\in \varphi(\tCHX)$ for any $l\geq0$, $r>0$. Furthermore,
		by \eqref{eq:haB-} and \eqref{eq:hax-}, we have
		$\gr(\haB_{[i,j],-l}*[K_{l\de-\alpha_{ij}}])\in \varphi(\tCHX)$ for any $l>0$. It follows that $\tCMHX^{\gr}=\varphi(\tCHX)\otimes\Q(\sqq)\big[[K_i]^{\pm1}\mid i\in\I\big]$ since $\tCMHX^{\gr}$ is generated by \eqref{eq:generators}. Together with the injectivity of $\varphi$,
		\eqref{eq:isoHiH} follows.
	\end{proof}

	\subsection{A presentation of $\tCMHX$}
	
	Recall that ${}^{\imath}\widetilde{\cc}(\bfk C_{p_i})$ is the composition subalgebra of $\iH(\bfk C_{p_i})$. The morphism $\widetilde{\psi}_{C_{p_i}}: \tUi_\sqq(\widehat{\mathfrak{sl}}_{p_i})\rightarrow \iH(\bfk C_{p_i})$ induces an isomorphism
	$\widetilde{\psi}_{C_{p_i}}: \tUi_\sqq(\widehat{\mathfrak{sl}}_{p_i})\rightarrow {}^{\imath}\widetilde{\cc}(\bfk C_{p_i})$. Inspired by Lemma \ref{lem:Hxm}, for $r>0$ we denote
	\begin{align}
		\widehat{H}_{[i,0],r}=\frac{[r]_{\sqq}}{r} \sum_{|\lambda|=r} \bn(\ell(\lambda)-1)\frac{[S_{i,0}^{(\lambda)}]}{\big|\Aut(S_{i,0}^{(\lambda)})\big|}-\delta_{r, ev} \sqq^{-\frac{r}{2}} \frac{[r/2]_{\sqq}}{r} [K_{\frac{r}{2}\de}]\in {}^{\imath}\widetilde{\ch}(\bfk C_{p_i}).
	\end{align}

	Recall the morphism $\Omega:  {}^{\text{Dr}}\tUi_{|v=\sqq}\rightarrow\tCMHX$ obtained in Theorem  \ref{thm:morphi}.
	For $1\leq i\leq \bt$, set
	\begin{align}
		\label{eq:Bs0}
		B_{[i,0]}:=\big[B_{[i,1],-1},B_{[i,2]},\cdots,B_{[i,p_i-1]}\big]_{v} C\K_{\theta_i}^{-1},
	\end{align}
	where $\theta_i:=\alpha_{i,1}+\cdots+\alpha_{i,p_i-1}$.

	\begin{lemma}
		We have
		$\Omega(B_{[i,0]})=\frac{-1}{q-1}[S_{i,0}]$ for any $1\leq i\leq \bt$. 
	\end{lemma}
	
	\begin{proof}
		A direct computation shows that
		\begin{align*}
			\big[[S_{i,0}^{(p_i-j)}],[S_{i,j+1}]\big]_{\sqq}=(1-q)[S_{i,0}^{(p_i-j-1)}]*[K_{S_{i,j+1}}]
		\end{align*}
		for any $1\leq j< p_i-1$.
		Then
		\begin{align*}
			\Omega(B_{[i,0]})=&(\frac{-1}{q-1})^{p_i-1}\big[[S_{i,0}^{(p_i-1)}],[S_{i,2}],\cdots, [S_{i,p_i-1}]\big]_{\sqq} * [K_{\theta_i-\alpha_{i,1}}]^{-1}
			\\
			=&(\frac{-1}{q-1})^{p_i-2}\big[[S_{i,0}^{(p_i-2)}], [S_{i,3}],\cdots, [S_{i,p_i-1}]\big]_{\sqq} * [K_{\theta_i-\alpha_{i,1}-\alpha_{i,2}}]^{-1}
			\\
			=&\cdots=\frac{-1}{q-1}[S_{i,0}].
		\end{align*}
	\end{proof}

	By using the presentation described in Corollary \ref{cor:presentationU+}, we have the following proposition.

	\begin{proposition}
		\label{prop:present i-comp}
		Let $\Gamma$ be of finite or affine type. Then $\tCMHX$ is isomorphic to the algebra ${}^\imath\widetilde{\bC}_{\X_\bfk}$ generated by ${}^{\imath}\widetilde{\ch}(\bfk C_{p_i})$, $[\co(l\vec{c})]$, $\widehat{H}_{\star,r}$, $[K_{\mu}]^{\pm1}$ for $1\leq i\leq \bt$,  $l\in\Z$, $r>0$,  $\mu\in\Z\I$, subject to the following relations precisely:
		\begin{align}
			\label{eq:HaDr0}
			&[K_\mu]*[K_{\mu}]^{-1}=1, \qquad [K_\mu] \text{ is central }, \text{ for }\mu\in\Z\I,
			\\
			\label{eq:HDr1}
			&\widehat{H}_{\star,r}=\widehat{H}_{[i,0],r}, \qquad\big[ {}^{\imath}\widetilde{\cc}(\bfk C_{p_i}), {}^{\imath}\widetilde{\cc}(\bfk C_{p_j}) \big]=0, \text{ for } 1\leq i\neq j\leq \bt,
			\\
			\label{eq:HaDr4}
			&\big[[\co(l\vec{c})],[S_{ij}]\big]=0, \text{ if }j> 1,
			\\
			\label{eq:HaDr2}
			&\big[\haH_{\star,m}, [\co(l\vec{c})]\big] =\frac{[2m]_\sqq}{m} [\co((l+m)\vec{c})]-\frac{[2m]_\sqq}{m} [\co((l-m)\vec{c})]*[K_{m\de}],
			\\
			\label{eq:HaDr3}
			&\big[[\co(k\vec{c})], [\co((l+1)\vec{c})]\big]_{\sqq^{-2}}  -\sqq^{-2} \big[[\co((k+1)\vec{c})], [\co(l\vec{c})]\big]_{\sqq^{2}}
			\\\notag
			&= (1-q)^2\Big(\sqq^{-2}\haTh_{\star,l-k+1} *[K_\de]^k* [K_{\co}]-\sqq^{-4} \haTh_{\star,l-k-1}* [K_\de]^{k+1} *[K_{\co}]
			\\\notag
			&+\sqq^{-2}\haTh_{\star,k-l+1}*[K_\de]^l *[K_{\co}]-\sqq^{-4} \haTh_{\star,k-l-1}*[ K_\de]^{l+1} *[K_\co]\Big),
			\\
			\label{eq:HaDr5}
			& \big[[S_{i,0}],[\co(l\vec{c})]\big]_{\sqq}=(1-q)^{2-p_i} \big[[\co((l-1)\vec{c})], [S_{i,1}],\cdots, [S_{i,p_i-1}]\big]_{\sqq}*[K_{S_{i,0}}],
			\\
			\label{eq:HaDr6}
			&\Big[[S_{i,0}], \big[[S_{i,1}],[\co(l\vec{c})]\big]_{\sqq^{-1}} \Big] =0, \text{ if } p_i>2,
			\\
			\label{eq:HaDr9}
			&\big[ [\co((l+1)\vec{c})],\haB_{[i,1],r}\big]_{\sqq^{-1}}-\sqq^{-1} \big[[\co(l\vec{c})],\haB_{[i,1],r+1}\big]_\sqq=0,\\
			\label{eq:HaDr7}
			&\widehat{\SS}(t_1,t_2|l; [i,1],\star) = (1-q)^2 \widehat{\R}(t_1,t_2|l; [i,1],\star)
			, \text{ for }t_1,t_2>0, l\in\Z,
			\\
			\label{eq:HaDr8}
			&\widehat{\SS}(l_1,l_2|r; [i,1],\star) = (1-q)^2 \widehat{\R}(l_1,l_2|r; \star,[i,1]),
			\text{ for } l_1,l_2\in\Z,  r>0.
		\end{align}
	\end{proposition}

	\begin{proof}
		We will show that the assignment
		$$[K_{\mu}]\mapsto [K_{\mu}],\, [\co(l\vec{c})]\mapsto [\co(l\vec{c})], \,[S_{ij}]\mapsto [S_{ij}], \, \widehat{H}_{\star,r} \mapsto \widehat{H}_{\star,r},\, \widehat{H}_{[i,0],r}\mapsto  \widehat{H}_{\star,r}$$ yields an algebra isomorphism
		\begin{align}
			\label{eq:isoUpsilon}
			&\Upsilon:{}^\imath\widetilde{\bC}_{\X_\bfk}\longrightarrow \tCMHX.
		\end{align}
		
		First, in order to prove that $\Upsilon$ is well defined, by using Corollary \ref{cor:present HCn} and the proof of Theorem \ref{thm:morphi}, it remains to check \eqref{eq:HaDr5}--\eqref{eq:HaDr6} in $\tCMHX$.
		
		A direct computation shows
		\begin{align}
			&\big[[S_{i,0}],[\co(l\vec{c})]\big]_{\sqq}=(1-q) [\co((l-1)\vec{c} +(p_i-1)\vec{x}_i )]*[K_{S_{i,0}}],
			\\
			\label{extension for A2}
			&\big[[\co(l\vec{c}+(j-1)\vec{x}_i)], [S_{ij}] \big]_{\sqq}=(1-q) [\co(l\vec{c} +j\vec{x}_i )],
		\end{align}
		for $l\in\mathbb{Z}$ and $1\leq j\leq p_i-1$.
		It follows that
		\begin{align*}
			[\co((l-1)\vec{c} +(p_i-1)\vec{x}_i )]=(1-q)^{-p_i+1} \big[[\co((l-1)\vec{c})],[S_{i,1}],\cdots, [S_{i,p_i-1}]\big]_{\sqq}
		\end{align*}
		Thus
		\begin{align*}
			\big[[S_{i,0}],[\co(l\vec{c})]\big]_{\sqq}=(1-q)^{2-p_i}
			\big[[\co((l-1)\vec{c})],[S_{i,1}],\cdots, [S_{i,p_i-1}]\big]_{\sqq}
			*[K_{S_{i,0}}],
		\end{align*}
		and then \eqref{eq:HaDr5} holds.
		
		Now we prove \eqref{eq:HaDr6}.
		By \eqref{extension for A2} we have
		\begin{align*}
			\big[[S_{i,1}],[\co(l\vec{c})]\big]_{\sqq^{-1}}=(\sqq-\sqq^{-1})[\co(l\vec{c}+\vec{x}_i)].
		\end{align*}
		Moreover, there are no homomorphisms and extensions between $\co(l\vec{c}+\vec{x}_i)$ and $S_{i,0}$ when $p_i>2$.
		Hence \eqref{eq:HaDr6} follows.
		
		Secondly, in order to prove that $\Upsilon$ is isomorphic, it is enough to prove that it is injective, which is given below.
		
		Define $\deg ([\co(l\vec{c})])=l\de+\alpha_\star$ for $l\in\Z$, $\deg (\widehat{H}_{\star,r})=r\de$ for $r>0$, and $\deg ([K_\mu])=0$ for $\mu\in K_0(\coh(\X_\bfk))$. Using the filtered algebra structure on ${}^{\imath}\widetilde{\ch}(\bfk C_{p_i})$, we can make ${}^\imath\widetilde{\bC}_{\X_\bfk}$ to be a filtered algebra. In this way, $\Upsilon$ is a morphism of  filtered  algebras. Denote by $\gr\Upsilon: {}^\imath\widetilde{\bC}_{\X_\bfk}^{\gr} \rightarrow \tCMHX^{\gr}$ the graded morphism between their associated graded algebras. 
		
		Using Lemma \ref{lem:phirestrction} and Corollary \ref{cor:presentationU+}, we have
		$$\rho': \tCMHX^{\gr}\cong \tCHX\otimes\Q(\sqq)\big[[K_i]^{\pm1}\mid i\in\I\big] \cong
		\widetilde{\bC}_{\X_\bfk}\otimes \Q(\sqq)\big[[K_i]^{\pm1}\mid i\in\I\big].$$
		By comparing the presentations of $\widetilde{\bC}_{\X_\bfk}$ and ${}^\imath\widetilde{\bC}_{\X_\bfk}$ obtained in Corollary \ref{cor:presentationU+} and Proposition \ref{prop:present i-comp} respectively, there exists a natural algebra morphism
		$$\rho'': \widetilde{\bC}_{\X_\bfk}\otimes \Q(\sqq)\big[[K_i]^{\pm1}\mid i\in\I\big] \longrightarrow {}^\imath\widetilde{\bC}_{\X_\bfk}^{\gr}.$$
		Then we have an algebra morphism
		$\rho:=\rho''\circ\rho': \tCMHX^{\gr}\rightarrow {}^\imath\widetilde{\bC}_{\X_\bfk}^{\gr}$.  Using \cite[Chapter 7, Corollary 6.14]{MR01}, we get that ${}^\imath\widetilde{\bC}_{\X_\bfk}^{\gr}$ is generated by ${}^{\imath}\widetilde{\ch}(\bfk C_{p_i})$, $[\co(l\vec{c})]$, $\widehat{H}_{\star,r}$, $[K_{\mu}]^{\pm1}$ for $1\leq i\leq \bt$,  $l\in\Z$, $r>0$,  $\mu\in\Z\I$.  Then it is easy to see that $\rho\circ\gr\Upsilon=\Id$. So $\gr\Upsilon$ is injective. Then $\Upsilon$ is injective by a standard argument of filtered algebras.
	\end{proof}
	
	\subsection{Injectivity of $\Omega$}
	In this subsection, we prove that $\Omega: \tUiD_{|v=\sqq}\rightarrow \tCMHX$ is an isomorphism if $\Gamma$ is of finite type or affine type.


	\begin{proposition}
		\label{prop:inverse}
		There exists an  algebra homomorphism $\Psi: {}^\imath\widetilde{\cc}(\X_\bfk)\rightarrow \tUiD_{|v=\sqq}$ which sends
		\begin{align}
			&[K_\mu]\mapsto \K_\mu, \quad\text{ for } \mu\in\Z\I,
			\\
			&[\co(l\vec{c})]\mapsto (1-q)B_{\star,l}, \qquad
			\haH_{\star,r}\mapsto H_{\star,r},\quad \text{ for }l\in\Z, r>0,
			\\
			&[S_{ij}]\mapsto (1-q)B_{[i,j]},\qquad [S_{i,0}]\mapsto (1-q)B_{[i,0]}, \quad \text{ for }1\leq i\leq \bt, 1\leq j\leq p_i-1.
		\end{align}
	\end{proposition}

	\begin{proof}
		By Proposition \ref{prop:present i-comp}, we identify $\tCMHX$ with ${}^\imath\widetilde{\bC}_{\X_\bfk}$ in the following. So it is enough to prove that $\Psi$ preserves the relations \eqref{eq:HaDr2}, \eqref{eq:HaDr3}, \eqref{eq:HaDr0}--\eqref{eq:HaDr8}. By Definition \ref{def:iDR},  \cite[(7.12)]{LR21} and \eqref{eq:Bs0} , it remains to prove that $\Psi$ preserves \eqref{eq:HaDr5}--\eqref{eq:HaDr6}, which follows from Lemma \ref{lem:relaitionS0} below.
	\end{proof}
	
	\begin{lemma}
		\label{lem:relaitionS0}
		For $1\leq i\leq \bt$ and $l\in\Z$, the following equalities hold in $\tUiD$:
		\begin{align}
			\label{eq:BBs0}
			&[ B_{[i,0]},B_{\star,l}]_{v}=-\big[ B_{\star,l-1},B_{[i,1]},B_{[i,2]},\cdots,B_{[i,p_i-1]}\big]_{v}C\K_{\theta_i}^{-1};  
			\\
			\label{eq:Bs0Bs1B}
			&\big[B_{[i,0]}, [B_{[i,1]},B_{\star,l}]_{v^{-1}}\big]=0,  \text{ if }p_i>2.
		\end{align}
	\end{lemma}
	
	\begin{proof}
		For \eqref{eq:BBs0}, by \eqref{iDR4}, we have
		$[B_{[i,j]},B_{\star,l}]=0$ for $2\leq j<p_i$. Then
		\begin{align*}
			\big[B_{[i,0]}, B_{\star,l} \big]_{v}=& \big[B_{[i,1],-1},B_{[i,2]}, \cdots,B_{[i,p_i-1]},B_{\star,l} \big]_{v}C\K_{\theta_i}^{-1}
			\\
			=& \big[[B_{[i,1],-1},B_{\star,l}]_{v}, B_{[i,2]},\cdots, B_{[i,p_i-1]}\big]_{v}C\K_{\theta_i}^{-1}
			\\
			=&- \big[B_{\star,l-1},B_{[i,1]}, B_{[i,2]},\cdots, B_{[i,p_i-1]}\big]_{v}C\K_{\theta_i}^{-1},
		\end{align*}
		where the third equality follows by  \eqref{iDR3a}.

		
		For \eqref{eq:Bs0Bs1B}, we have
		\begin{align*}
			\big[B_{[i,0]}, [B_{[i,1]},B_{\star,l}]_{v^{-1}}\big]
			=&-\big[B_{\star,l},[B_{i,0},B_{i,1}]_{v^{-1}}\big]-\big[B_{[i,1]},[B_{\star,l},B_{[i,0]}]_{v^{-1}}\big].
		\end{align*}
		Claim:
		\begin{align}
			\label{claim:BBB1}
			\big[B_{\star,l},[B_{i,0},B_{i,1}]_{v^{-1}}\big]=0,
			\\
			\label{claim:BBB2}
			\big[B_{[i,1]},[B_{\star,l},B_{[i,0]}]_{v^{-1}}\big]=0.
		\end{align} Then \eqref{eq:Bs0Bs1B} follows.

		For \eqref{claim:BBB1}, since $p_i>2$, by \eqref{iDR4} and \eqref{iDR3a}, we have
		\begin{align}
			\label{eq:reduction}
			[B_{[i,0]},B_{[i,1]}]_{v^{-1}}=&\Big[  \big[B_{[i,1],-1},B_{[i,2]},\cdots,B_{[i,p_i-1]}\big]_{v} ,B_{[i,1]}\Big]_{v^{-1}}C\K_{\theta_i}^{-1}
			\\\notag
			=&\Big[  \big[[B_{[i,1],-1},B_{[i,2]}]_v,B_{[i,1]}\big]_{v^{-1}},\cdots,B_{[i,p_i-1]}\Big]_{v} C\K_{\theta_i}^{-1}
			\\\notag
			=&-\Big[  \big[[B_{[i,2],-1},B_{[i,1]}]_v,B_{[i,1]}\big]_{v^{-1}},\cdots,B_{[i,p_i-1]}\Big]_{v} C\K_{\theta_i}^{-1}
			\\\notag
			=&v^{-1}\Big[B_{[i,2],-1},B_{[i,3]},\cdots,B_{[i,p_i-1]}\Big]_{v} \K_{[i,1]} C\K_{\theta_i}^{-1},
		\end{align}
		where we use \eqref{iDR5} in the last equality.
		Using \eqref{iDR4} again, we have $\big[B_{\star,l},[B_{i,0},B_{i,1}]_{v^{-1}}\big]=0$.

		For \eqref{claim:BBB2}, using \eqref{eq:BBs0}, we have
		\begin{align*}
			&\big[B_{[i,1]},[B_{\star,l},B_{[i,0]}]_{v^{-1}}\big]
			=v^{-1}\Big[B_{[i,1]}, \big[B_{\star,l-1},B_{[i,1]}, B_{[i,2]},\cdots, B_{[i,p_i-1]}\big]_{v} \Big]C\K_{\theta_i}^{-1}.
		\end{align*}
		Using $\big[ B_{[i,1]}, B_{[i,j]}\big]=0$ for $j>2$, it is enough to show that
		\begin{align}
			\label{eq:BBBB}
			\Big[\big[B_{\star,l-1}, B_{[i,1]},B_{[i,2]}\big]_{v},B_{[i,1]}\Big]=0.
		\end{align}
		In fact, since $[B_{\star,l-1},B_{i,2}]=0$, we have
		\begin{align*}
			&(-v^{-2})\cdot\text{LHS}\eqref{eq:BBBB}\\
			= & B_{[i,1]} B_{[i,2]} B_{[i,1]}B_{\star,l-1} -v^{-1} B_{[i,1]} B_{[i,1]} B_{\star,l-1} B_{[i,2]} +v^{-2}B_{[i,1]}B_{\star,l-1} B_{[i,1]}B_{[i,2]}
			\\
			&-B_{[i,2]} B_{[i,1]}B_{\star,l-1}B_{[i,1]} +v^{-1}B_{[i,2]}B_{\star,l-1}B_{[i,1]}B_{[i,1]} -v^{-2} B_{\star,l-1}B_{[i,1]}B_{[i,2]}B_{[i,1]}
			\\
			=&\frac{1}{[2]} \big( B_{[i,1]}  B_{[i,1]}B_{[i,2]} +B_{[i,2]}B_{[i,1]}  B_{[i,1]}+v^{-1}B_{[i,2]}\K_{[i,1]}\big)B_{\star,l-1}-v^{-1} B_{[i,1]} B_{[i,1]}  B_{[i,2]}B_{\star,l-1}
			\\
			&+\frac{v^{-2}}{[2]}\big( B_{[i,1]}B_{[i,1]}B_{\star,l-1}+ B_{\star,l-1}B_{[i,1]}B_{[i,1]}+v^{-1}B_{\star,l-1}\K_{[i,1]}\big)B_{[i,2]}
			\\
			&-\frac{1}{[2]} B_{[i,2]} \big( B_{[i,1]}B_{[i,1]}B_{\star,l-1}+ B_{\star,l-1}B_{[i,1]}B_{[i,1]}+v^{-1}B_{\star,l-1}\K_{[i,1]}\big)+v^{-1}B_{[i,2]}B_{\star,l-1}B_{[i,1]}B_{[i,1]}
			\\
			&-\frac{v^{-2}}{[2]} B_{\star,l-1}\big( B_{[i,1]}  B_{[i,1]}B_{[i,2]} +B_{[i,2]}B_{[i,1]}  B_{[i,1]}+v^{-1}B_{[i,2]}\K_{[i,1]}\big)
			\\
			=&0.
		\end{align*}
		Here the second equality follows by using \eqref{iDR5}. The proof is completed.
	\end{proof}

	Now we can prove the main result Theorem \ref{thm:morphi}.

	\begin{proof}[Proof of Theorem \ref{thm:morphi}]
		It is enough to prove that 
		$\Omega: \tUiD_{|v=\sqq}\rightarrow \tCMHX$ is an isomorphism.  By definition, $\Omega$ is epimorphism. 	
		By Lemma \ref{lem:reduced generators} and Proposition \ref{prop:inverse}, it is easy to see that $\Psi\circ\Omega=\Id$, and then $\Omega$ is injective. So $\Omega: \tUiD_{|v=\sqq}\rightarrow \tCMHX$ is an isomorphism.
	\end{proof}

	\begin{remark}
		For the injectivity of $\Omega$, the corresponding result for the whole quantum loop algebras has not been proved except for the special case of the quantum affine $\mathfrak{sl}_2$. In fact, it is only proved in \cite{Sch04} that a ``positive part'' of the quantum loop algebras can be realized by the Hall algebra of $\X_\bfk$ if $\fg$ is of finite or affine type; see Theorem \ref{thm:Sch}.
		
		We expect the algebra epimorphism $\Omega:\tUiD_{|v=\sqq}\rightarrow \tCMHX$ to be injective for any Kac-Moody algebra $\fg$; see a similar expectation for quantum loop algebras in \cite{Sch04,DJX12}. 
	\end{remark}



\end{document}